\documentclass{preprint}

\ifoot{~}
\ofoot{L. Failer, T. Richter}
\ohead{\thepage}
\ihead{A parallel Newton-multigrid framework for monolithic
  fluid-structure interactions} 

\usepackage{graphicx}
\usepackage[table]{xcolor}
\usepackage{booktabs,units,xfrac}
\usepackage{placeins}
\usepackage{pgfplots}
\pgfplotsset{compat=1.9}
\usetikzlibrary{
  pgfplots.groupplots,
  matrix
}
\usepackage{amsmath,amssymb,dsfont,graphicx}

\usepackage[table]{xcolor}

\usepackage{algpseudocode}
\usepackage{algorithm}
\usepackage{float}
\floatstyle{plaintop}
\newfloat{algorithm}{tbp}{loa}[section]
\floatname{algorithm}{Algorithm}

\usepackage{tikz}
\usepackage{pgfplots}
\newlength \figureheight
\newlength \figurewidth

\newcommand{\tr}{\operatorname{tr}}

\newcommand{\vt}{\mathbf{v}}
\newcommand{\pt}{\mathbf{p}}

\newcommand{\ut}{\mathbf{u}}
\newcommand{\bt}{\mathbf{b}}
\newcommand{\Ft}{\mathbf{F}}
\newcommand{\Et}{\mathbf{E}}
\newcommand{\onet}{\mathbf{1}}

\newcommand{\deltat}{\mathbf{\delta}}

\newcommand{\sigmat}{\boldsymbol{\sigma}}
\newcommand{\Sigmat}{\boldsymbol{\Sigma}}

\newcommand{\ft}{\mathbf{f}}

\makeatletter
\newcommand{\opnorm}{\@ifstar\@opnorms\@opnorm}
\newcommand{\@opnorms}[1]{%
  \left|\mkern-1.5mu\left|\mkern-1.5mu\left|
   #1
  \right|\mkern-1.5mu\right|\mkern-1.5mu\right|
}
\newcommand{\@opnorm}[2][]{%
  \mathopen{#1|\mkern-1.5mu#1|\mkern-1.5mu#1|}
  #2
  \mathclose{#1|\mkern-1.5mu#1|\mkern-1.5mu#1|}
}
\makeatother

\newcommand{\SO}{\mathcal{S}}
\newcommand{\IN}{\mathcal{I}}
\newcommand{\At}{\mathbf{A}}
\newcommand{\FL}{\mathcal{F}}
\newcommand{\Ut}{\mathbf{U}}
\newcommand{\xt}{\mathbf{x}}

\newcommand{\norm}[1]{\lVert #1 \rVert}

\newcommand{\SCC}{\cellcolor{gray!20}}

\definecolor{mycolor1}{rgb}{0.00000,0.44700,0.74100}%
\definecolor{mycolor2}{rgb}{0.85000,0.32500,0.09800}%
\definecolor{mycolor3}{rgb}{0.92900,0.69400,0.12500}%
\definecolor{mycolor4}{RGB}{152,198,234}%
\definecolor{mycolor5}{RGB}{1,1,1}%
\begin{document}

\title{A parallel Newton multigrid framework for monolithic
  fluid-structure interactions
}

\author{L. Failer\thanks{
    Technische Universit\"at M\"unchen, 
    85748 Garching bei M\"unchen, Germany,
    \texttt{lukas.failer@ma.tum.de}
  }
  \and
  T. Richter\thanks{
    Otto-von-Guericke Universit\"at Magdeburg, 
    39104 Magdeburg, Germany, 
    \texttt{thomas.richter@ovgu.de}, 
    and Interdisciplinary Center for Scientific Computing, 
    Heidelberg University, 69120 Heidelberg, Germany
  }
}

\maketitle

\begin{abstract}
  We present a monolithic parallel Newton-multigrid solver for
nonlinear three dimensional fluid-structure interactions  in
Arbitrary Lagrangian Eulerian (ALE) formulation. We start with a finite
element discretization of the coupled problem, based on a remapping
of the Navier-Stokes equation onto a fixed reference framework. The
strongly coupled fluid-structure interaction problem is discretized
with finite elements in space and finite differences in time.
The resulting nonlinear and linear systems of equations are large
and show a very high condition number.

We present a novel Newton approach that is based on two essential
ideas: First, a condensation of the solid deformation by
exploiting the discretized  
velocity-deformation relation $d_t \ut = \vt$. Second, the Jacobian
of the fluid-structure interaction system is simplified by
neglecting all derivatives with respect to the ALE deformation, an
approximation that has shown to have little impact. The resulting
system of equations decouples into a joint momentum equation and into
two separated equations for the deformation fields in solid and
fluid. Besides a reduction of the problem sizes, the approximation
has a positive effect on the conditioning of the systems such that
multigrid solvers with simple smoothers like a parallel
Vanka-iteration can be applied.

We demonstrate the efficiency of the resulting solver infrastructure
on a well-studied 2d test-case and we also introduce a 
challenging 3d problem. For 3d problems we achieve a substantial
accelaration as compared to established approaches found in literature.

\end{abstract}

\section{Introduction}

Fluid structure interactions appear in various problems ranging from
classical applications in engineering like the design of ships or
aircrafts, the design of wind turbines, but they are also present in
bio/medical systems describing the blood flow in the heart or in
general problems involving the cardiovascular system. The typical
challenge of fluid-structure interactions is two-fold. First, the
special coupling character that stems from the coupling of a
hyperbolic-type equation - the solid problem - with a parabolic-type
equation - the Navier-Stokes equations. Second, the moving domain
character brings along severe nonlinearities that have a non-local
character, as geometrical changes close to the moving fluid-solid
interface might have big impact on the overall solution.

Numerical approaches can usually be classified into \emph{monolithic
  approaches}, where the coupled fluid-structure interaction system is
taken as one entity and into \emph{partitioned approaches}, where two
separate problems - for fluid and solid - are formulated and where the
coupling between them is incorporated in terms of an outer (iterative)
algorithm. This second approach has the advantage that difficulties
are isolated and that perfectly suited numerical schemes can be used
for each of the subproblems. There are however application classes
where partitioned approaches either fail or lack efficiency. The
\emph{added mass effect}~\cite{CausinGerbeauNobile2005} exactly
describes this special 
stiffness connected to fluid-structure interactions. It is typical for
problems with similar densities in the fluid and the solid - as it
happens in the interaction of blood and tissue or in the interaction
of water and the solid structure of a vessel. Here, monolithic
approaches are considered to be favourable.

Monolithic approaches all give rise to strongly coupled, usually very
large and nonlinear algebraic systems of equations. Although there has
been substantial progress in designing efficient numerical schemes for
tackling the nonlinear problems~\cite{HronTurek2006a,HeilHazelBoyle2008,FernandezGerbeau2009} (usually by Newton's method)
and the resulting linear systems~\cite{GeeKuettlerWall2010,TurekHronMadlikRazzaqWobkerAcker2010,Richter2015,LangerYang2017,AulisaBnaBornia2018,CrosettoDeparis,DeparisForti}, the computational effort is
still immense and numerically accurate results for 3d problems are
still rare.

In this contribution we present an approximated  Newton scheme for
solving nonstationary fluids structure interactions in a strictly
monolithic formulation. The idea is based on the observation that
the Newton convergence rate does not significantly worsen, if we
neglect the derivatives with respect to the ALE deformation,
see~\cite[Section 5.2.3]{Richter2017}. Although convergence rates
slightly suffer,  overall computational times can be reduced due to
lesser effort for assembling the matrix. Here, we exploit this
structure of the reduced Jacobian to achieve an exact splitting of the
monolithic Jacobian into a coupled problem for the velocities of fluid
and solid and into a second step, where separate update problems are
solved for  solid and fluid deformation. Apart from the approximation
of the Jacobian, no further splitting error is introduced. The benefit
of this approach is twofold: instead of one large system with 7
coupled unknowns (pressure, velocity field and deformation field in
3d) we solve one coupled system of four unknowns (pressure and
velocities) and two separate problems involving the
deformations of each domain. Second, separating a reduced velocity
problem has a 
positive effect on the system matrices such that efficient
preconditioners and smoothers can be applied that are suitable for
easy parallelization. Finally, we use the newly developed solver to
introduce and test a new three dimensional benchmark configuration
that is based on the configurations described by Hron and
Turek~\cite{HronTurek2006a}.

In the following section we give a brief presentation of the
fluid-structure interaction problem in a variational Arbitrary
Lagrangian Eulerian formulation. Section~\ref{sec:disc} shortly
presents the discretization of the equations in space and time. As
formulation and discretization are based on established techniques,
these two sections are rather concise. The nonlinear and linear
solution framework is described in Section~\ref{sec:solver}, where we
start by an approximation of the Jacobian that results in a
natural partitioning of the linear systems, which in turn are
approximated by parallel multigrid methods. Numerical test-cases
demonstrate the efficiency and scalability in
Section~\ref{sec:num}. Here, we also present a new and challenging
3d configuration for benchmarking fluid-structure interactions. We
conclude in Section~\ref{sec:conclusion}.

\section{Governing equations}

Here, we present the monolithic formulation for fluid structure
interactions, coupling the incompressible Navier-Stokes equations and
an hyperelastic solid, based on the St. Venant Kirchhoff material. For
details we refer to~\cite{Richter2017}. 

On the $d$-dimensional domain, partitioned in reference configuration
$\Omega = \FL\cup\IN\cup\SO$, where $\FL$ is the fluid domain, $\SO$
the solid domain and $\IN$ the fluid structure interface,  we denote
by $\vt$ the velocity field, split into fluid  velocity
$\vt_f:=\vt|_{\FL}$ and solid velocity $\vt_s:=\vt|_{\SO}$, and 
by $\ut$ the deformation field, again with $\ut_s:=\ut|_{\SO}$ and
$\ut_f:=\ut|_{\FL}$. The boundary of the fluid domain
$\Gamma_f:=\partial\FL\setminus\IN$ is split into inflow boundary
$\Gamma_f^{in}$ and wall boundary $\Gamma_f^{wall}$, where we usually
assume Dirichlet conditions,
$\Gamma_f^D:=\Gamma_f^{in}\cup\Gamma_f^{wall}$,  and a possible outflow boundary
$\Gamma_f^{out}$, where we enforce the do-nothing outflow
condition~\cite{HeywoodRannacherTurek1992}. The solid boundary
$\Gamma_s=\partial\SO\setminus\IN$ is split into Dirichlet part
$\Gamma_s^D$ and a Neumann part $\Gamma_s^N$. 

We formulate the coupled fluid-structure interaction problem in a
strictly monolithic scheme by mapping the moving fluid domain onto the
reference state via the ALE map $T_f(t):\FL\to\FL(t)$,
constructed by a fluid domain deformation $T_f(t)=\operatorname{id} +
\ut_f(t)$. In the solid domain, this map
$T_s(t)=\operatorname{id}+\ut_s(t)$ denotes the Lagrange-Euler
mapping and as the deformation field $\ut$ will be defined globally on
$\Omega$ we simply use the notation $T(t)=\operatorname{id}+\ut(t)$
with the deformation gradient $\Ft:=\nabla T$ and its determinant
$J:=\operatorname{det}(\Ft)$. 
We find the global (in fluid and solid domain)
velocity and deformation fields $\vt$ and $\ut$ and the  pressure $p$
in the function 
spaces 
\[
\vt(t)\in \vt^D(t)+H^1_0(\Omega;\Gamma_f^D\cup\Gamma_s^D)^d,\quad
\ut(t)\in
\ut^D(t)+H^1_0(\Omega;(\partial\FL\setminus\IN)\cup\Gamma_s^D)^d,\quad
p\in L^2(\FL)
\]
as solution to 
\begin{equation}\label{aletime:1}
  \begin{aligned}
    \big( J(\partial_t  \vt +  (\Ft^{-1}( \vt-\partial_t   
    \ut)\cdot\nabla) \vt,\phi\big)_{\FL}
    + \big( J\sigmat_f
    \Ft^{-T},\nabla\phi\big)_{\FL}\hspace{-1cm} \\
   +(\rho_s^0\partial_t \vt,\phi)_{\SO}
    +( \Ft\Sigmat_s,\nabla\phi)_{\SO}
    &=
    ( J\rho_f \ft,\phi)_{\FL} 
    +(\rho_s^0 \ft,\phi)_{\SO}\\
    \big(J\Ft^{-1}:\nabla\vt^T,\xi\big)_{\FL} &= 0\\
    (\partial_t \ut- \vt,\psi_s)_{\SO}&=0\\
    (\nabla \ut,\nabla\psi_f)_{\FL}&=0,
  \end{aligned}
\end{equation}
where the test functions are given in 
\[
\phi\in H^1_0(\Omega;\Gamma_f^D\cup\Gamma_s^D)^d,\quad
\xi\in L^2(\FL),\quad
\psi_f\in H^1_0(\FL)^d,\quad
\psi_s\in L^2(\SO)^d.
\]
By $\rho_s^0$ we denote the solid's density, by $\ut^D(t)\in
H^1(\Omega)^d$ and $\vt^D(t)\in H^1(\Omega)^d$  extensions of the
Dirichlet data into the domain.  
The Cauchy stress
tensor of the Navier-Stokes equations in ALE coordinates is given by
\[
\sigmat_f(\vt,p) = -p_f I + \rho_f\nu_f (\nabla\vt\Ft^{-1} + \Ft^{-T}\nabla\vt^T)
\]
with the kinematic viscosity $\nu_f$ and the density $\rho_f$. In the
solid we consider the St. Venant Kirchhoff 
material with the Piola Kirchhoff tensor
\[
\Sigmat_s(\ut) = 2\mu_s \Et_s + \lambda_s\tr(\Et_s)I,\quad
\Et_s:=\frac{1}{2}(\Ft^T\Ft-I)
\]
and with the shear modulus $\mu_s$ and the Lam\'e coefficient
$\lambda_s$. In~(\ref{aletime:1}) we construct the ALE extension
$\ut_f=\ut|_\FL$ by a simple harmonic extension. A detailed discussion
and further literature on the construction of this extension is found
in~\cite{YirgitSchaeferHeck2008,Richter2017}.

For shorter notation, we denote by $U:=(\vt,\ut,p_f)$ the solution
and by $\Phi:=(\xi,\phi,\psi_f,\psi_s)$ the test functions.

\section{Discretization}\label{sec:disc}

We give a very brief presentation on the numerical approximation of
System~(\ref{aletime:1}). In time, we use the theta time stepping
scheme, which includes the backward Euler method, the Crank-Nicolson
scheme and variants like the fractional step theta method,
see~\cite{TurekRivkindHronGlowinski2006}. In space we use conforming
finite elements. 

\subsection{Temporal discretization}

For discretization in time we split the temporal interval $I=[0,T]$
into discrete time steps $0=t_1<t_2<\cdots < t_N = T$ with the step
size $k:=t_n-t_{n-1}$. For simplicity we assume that the subdivision
is uniform. By $U_n\approx U(t_n)$ we denote the approximation at time
$t_n$. We choose the theta time stepping method for temporal
discretization with $\theta\in [0,1]$.
To simplify the presentation we introduce
\begin{equation}\label{aletime:1.5}
  \begin{aligned}
    A_F(U,\phi) &:= \big(J(\Ft^{-1}\vt\cdot\nabla)\vt,\phi\big)_\FL
    + \big(\rho_f\nu_f
    J(\nabla\vt\Ft^{-1}+\Ft^{-T}\nabla\vt^T)\Ft^{-T},\nabla\phi\big)_\FL
    -   \big(J\rho_f\ft,\phi\big)_\FL \\
    A_S(U,\phi) &:=  \big(\Ft\Sigmat_s,\nabla\phi\big)_{\SO}
    -\big(\rho_s^0\ft,\phi\big)_\SO ,\quad
    A_{ALE}(U,\psi_f) := \big(\nabla\ut,\nabla\psi_f\big)_{\FL}    
    \\
    A_{p}(U,\phi) &:=  \big(Jp\Ft^{-1},\nabla\phi\big)_{\FL} ,\quad
    A_{div}(U,\xi) := \big( J\Ft^{-1}:\nabla\vt^T,\xi\big)_{\FL}. 
  \end{aligned}
\end{equation}
Then,  one time step $t_{n-1}\mapsto t_n$ of the theta scheme is given as 
\begin{equation}\label{aletime:2}
  \begin{aligned}
    \underbrace{ \big(\bar J_n (\vt_n-\vt_{n-1}),\phi\big)_{\FL}
    - \big(( \bar J_n\bar\Ft^{-1}
    (\ut_n-\ut_{n-1})\cdot\nabla)\bar\vt_n,\phi\big)_{\FL}}_{F_{NS}(U_n,\phi)}&
    + \underbrace{k A_p(U_n,\phi) +k\theta A_F(U_n,\phi) }_{F_{NS}(U_n,\phi)}\quad \\
    +\big(\rho^0_s (\vt_n-\vt_{n-1}),\phi\big)_{\SO}+k\theta A_S(U_n,\phi)
    &= - k(1-\theta) A_F(U_{n-1},\phi) \\
    &\qquad - k(1-\theta) A_S(U_{n-1},\phi) \\
    k A_{div}(U_n,\xi) &=0\\
    k A_{ALE}(U_n,\psi_f)& =0\\
    \big(\ut_n,\psi_s\big)_{\SO}
    -k\theta\big(\vt_n,\psi_s\big)_\SO =
    \big(\ut_{n-1},\psi_s\big) + &
    k(1-\theta)\big(\vt_{n-1},\psi_s\big)_\SO,\\
  \end{aligned}
\end{equation}
with $\bar J_n = \sfrac{1}{2}(J_{n-1}+J_n)$ and $\bar\Ft_n =
\sfrac{1}{2}(\Ft_{n-1} + \Ft_n)$. 
Note that the ALE extension equation $A_{ALE}$, the divergence
equation $A_{div}$ and the
pressure coupling $A_p$ are completely implicit. A discussion of this
scheme and results on its stability for fluid-structure interactions
are found in~\cite{RichterWick2015_time,Richter2017}. Usually we
consider $\theta = \sfrac{1}{2}+{\cal O}(k)$ to get second order
convergence and good stability properties.

The last equation in~(\ref{aletime:2}) gives a relation for the new
deformation at time $t_n$
\[
\ut_n = \ut_{n-1}+k\theta \vt_n + k(1-\theta)\vt_{n-1}\text{ in }\SO
\]
and we will use this representation to eliminate the unknown
deformation and base the 
solid stresses purely on last time step and the unknown velocity, i.e.
by expressing the deformation gradient as 
\begin{equation}\label{defgrad}
  \Ft_n=\Ft(\ut_{n}) \,\widehat{=}\, \Ft(\ut_{n-1},\vt_{n-1};\vt_n) =
  I+\nabla   \big(\ut_{n-1}+k\theta \vt_n + 
  k(1-\theta)\vt_{n-1}\big)\text{ in }\SO. 
\end{equation}
Removing the solid deformation from the momentum equation will help to
reduce the algebraic systems in Section~\ref{sec:solver}. A similar
technique within a Eulerian formulation and using a characteristics
method is presented in~\cite{Pironneau2016,Pironneau2019}.

\subsection{Finite elements}

In space, we discretize with conforming finite elements by choosing
discrete function spaces $U_h\in X_h$ and $\Phi_h\in Y_h$.
We only consider finite element meshes that resolve the
interface $\IN$ in the reference configuration, such that the ALE
formulation will always exactly track the moving interface. 
In our
setting, implemented in the finite element library Gascoigne
3D~\cite{Gascoigne3D} we use quadratic finite elements for all
unknowns and add stabilization terms based on local
projections~\cite{BeckerBraack2001,Frei2016,Molnar2015,Richter2017} to
satisfy the inf-sup condition. Where transport is dominant,
additional stabilization terms of streamline upwind
type~\cite{Wall1999,RichterWick2010,HronTurek2006a} or of local projection
type~\cite{Richter2017,Failer2017}  are added. As the remainder of
this manuscript only considers the fully discrete setting, we refrain
from indicating spatial or temporal discrete variables with the usual
subscripts. 

For each time step $t_{n-1}\mapsto t_n$ we introduce the following
short notation for the system of algebraic equations that is based on
the splitting of the solution into unknowns acting in the fluid domain
$(\vt_f,\ut_f)$, on the interface $(\vt_i,\ut_i)$  and those on the solid
$(\vt_{s},\ut_{s})$.  The pressure variable $p$ acts in the fluid
and on the interface.
\begin{equation}\label{system}
  {\cal A}(U): = \begin{pmatrix}
    {\cal D}(p,\vt_f,\ut_f,\vt_{i},\ut_{i},\vt_{s},\ut_{s}) \\    
    {\cal M}^f(p,\vt_f,\ut_f,\vt_{i},\ut_i) \\
    {\cal M}^i(p,\vt_f,\ut_f,\vt_{i},\ut_i,\vt_s) \\
    {\cal M}^s(p,\vt_{i},\ut_i,\vt_s) \\
    {\cal E}(\ut_f,\ut_i)\\
    {\cal U}^i(\vt_{i},\ut_{i},\vt_s,\ut_s)\\
    {\cal U}^s(\vt_{i},\ut_{i},\vt_s,\ut_s)
  \end{pmatrix}
  =\begin{pmatrix}
  {\cal B}_1 \\ {\cal B}_2\\ {\cal B}_3\\ {\cal B}_4\\
  {\cal B}_5\\ {\cal B}_6 \\{\cal B}_7
  \end{pmatrix}=: {\cal B}
\end{equation}
${\cal D}$ describes the divergence equation which acts in the
fluid domain and on the interface, ${\cal M}$
the two momentum equations, acting in the fluid domain, on the
interface and in the solid domain (which is indicated by a
corresponding index), ${\cal E}$ describes the ALE extension in the fluid domain
and ${\cal U}$ is the relation between solid velocity and
solid deformation, acting on the interface degrees of freedom and in
the solid. Note that ${\cal M}^i$ and ${\cal M}^s$, the term describing the
momentum equations, do not directly depend on the solid deformation 
$\ut_{s}$ as we base the deformation gradient on the velocity,
see~(\ref{defgrad}).

\section{Solution of the algebraic systems}\label{sec:solver}

In fluid-structure interactions the solid and fluid problem are
coupled via interface conditions. Forces in normal direction along the
interface have to be equal (dynamic coupling condition) and the fluid
domain has to follow the solid motion (kinematic and geometric
coupling condition). If the solid motion is rather small and slow the
energy exchange happens mainly via the dynamic coupling
conditions. This allows the use of explicit time-stepping schemes for
the mesh motion and ALE transformation for these examples. We want to
follow a different approach and use a fully implicit time stepping
with an inexact Jacobian in the Newton algorithm. We neglect the
derivatives with respect to the ALE deformation. Thereby, we have to
solve in every Newton step a linear system of the same complexity as
in the case of a partitioned  time-stepping scheme. 

In~\cite[chapter 5]{Richter2017} we give a numerical study on
different linearization techniques. It is found that the overall
computational time can be reduced by neglecting the ALE derivatives in
the Jacobian. Even for the fsi-3 benchmark problem of Hron and
Turek~\cite{HronTurek2006} it is more efficient (in terms of overall
computational time) to omit these derivatives at the cost of some
additional Newton steps. Neglecting the ALE derivatives will be
crucial for the reduction step described in the following section.

As we only change the Jacobian, we still apply a fully implicit
time-stepping scheme and take advantage of its stability
properties. Furthermore the transport due to the mesh motion is well
approximated. For  small time-step sizes we will still observe
super-linear convergence as with an exact Newton algorithm. In
addition, the simplified structure of the matrix simplifies the
development of preconditioners sincerely as we will see later.

\subsection{Relation to approaches in literature}

Many (perhaps most) works on solvers for fluid-structure interactions
are based on partitioned schemes, where highly tuned schemes can be
applied to the two subproblems and acceleration methods are developed
for the coupling. For an overview on some methods we refer to
contributions in~\cite{BungartzSchaefer2006,BungartzSchaefer2010} and
the literature cited therein. We focus on problems with a dominant
added mass effect, where monolithic approaches are believed to be more
efficient~~\cite{HeilHazelBoyle2008}.

In the following we assume that the monolithic problem is approximated
with a Newton scheme. It has been
documented~\cite{Richter2015,AulisaBnaBornia2018} that the Jacobian is
very ill-conditioned with condition numbers exceeding those in fluid
or solid mechanics by far. Furthermore, the systems are (in particular
in 3d) so large that direct solvers are not applicable. In addition we 
found~\cite{Richter2015} that the condition numbers may be so large
that direct solvers do not even converge well.\footnote{These results
  where 
found in~\cite{Richter2015} for the direct solver
UMFPACK~\cite{Davis2014}. As similar study
in~\cite{AulisaBnaBornia2018} could validate our estimates for the
condition numbers but found better performance in the solver 
MUMPS~\cite{MUMPS:2}.}
All successfull solution strategies will
therefore feature some kind of partitioning, usually be means of a
decoupled preconditioner within a GMRES
iteration. In~\cite{LangerYang2017} an overview on state of the art
precondition techniques for iterative fluid-structure interaction
solvers is given. 

Multigrid solvers have first been used to accelerate the solution of
the subproblems within an iterative scheme. A fully monolithic
geometric multigrid
approach was presented in~\cite{HronTurek2006a} for 2d fsi
problems. Here, the multigrid smoother was based on a Vanka
iteration.
In~\cite{BrummelenZeeBorst2008} the authors analyzed a highly
simplified model problem and showed that a partitioned iteration as
smoother should result in ideal multigrid performance with improved
convergence rates on deeper mesh hierarchies. 
An algebraic multigrid method with applications in 2d and
3d~\cite{GeeKuettlerWall2010} was based on a Gauss-Seidel splitting in
the smoother. In~\cite{Richter2015} we presented a fully geometric
monolithic multigrid method with a smoother that is based on a
partitioning into fluid and solid problem and a block decomposition of
each equation. This approach has been extended to
incompressible materials and also to direct-to-steady-state
solutions~\cite{AulisaBnaBornia2018}. 

Some of these contributions employ parallelism. Recently, a
block-preconditioned parallel GMRES iteration was
presented~\cite{JodlbauerLanger}  and showed good performance on
various 2d and 3d test cases. A Gauss-Seidel decoupling with highly
efficient and massively parallel preconditioners based on the SIMPLE
scheme for the fluid and multigrid for a linear elasticity problem is
presented in~\cite{DeparisForti}.

\subsection{Linearization and splitting}

Each time step of the fully discrete problem is solved by Newton's
method. Evaluating the Jacobian is cumbersome due to the moving domain
character of the fluid problem. First presentations of the derivatives
of the fsi problem with respect to the mesh motion based on the
concept of shape derivatives have been given by Fernandez and
Moubachir~\cite{FernandezMoubachir2005}. Details in the spirit of our
formulation in ALE coordinates are given in~\cite[Section
  5.2.2]{Richter2017}.  
Based on the notation~(\ref{system}) let $U^{(0)}$ be an
initial guess (usually taken from the last time step) we iterate for
$l=0,1,2,\dots$
\begin{equation}\label{newton}
  {\cal A}'(U^{(l)}) W^{(l)} = {\cal B}-{\cal A}(U^{(l)}),\quad
  U^{(l+1)}:=U^{(l)}+\omega^{(l)} \cdot W^{(l)},
\end{equation}
with a line search parameter $\omega^{(l)}>0$ and the Jacobian ${\cal
  A}'(U)$ evaluated at $U$. Each linear problem can be written as 
\begin{equation}\label{FULLMATRIX}
  \left(
  \begin{array}{ccc|cc|cc}
    0  & {\cal D}_{\vt_f} & \SCC {\cal D}_{\ut_f} & {\cal D}_{\vt_i} & \SCC {\cal D}_{\ut_i} & 0 & 0\\
    {\cal M}^f_p & {\cal M}^f_{\vt_f} & \SCC {\cal M}^f_{\ut_f} & {\cal M}^f_{\vt_i} & \SCC {\cal M}^f_{\ut_i}& 0 & 0 \\
    {\cal M}^i_p & {\cal M}^i_{\vt_f} & \SCC {\cal M}^i_{\ut_f} &
    {\cal M}^i_{\vt_i} &\boldsymbol{ \SCC {\cal M}^i_{\ut_i}}& {\cal
      M}^i_{\vt_s} & \boldsymbol{{\cal M}^i_{\ut_s}}\\
    {\cal M}^s_p & 0  & 0 & {\cal M}^s_{\vt_i} & \boldsymbol{ {\cal
        M}^s_{\ut_i}}& {\cal M}^s_{\vt_s} & \boldsymbol{{\cal M}^s_{\ut_s}}\\
    \hline
    0 & 0 & {\cal E}^f_{\ut_{f}} & 0 & {\cal E}^f_{\ut_{i}}&0&0\\
    \hline
    0&0&0&{\cal U}^i_{\vt_{i}}&{\cal U}^i_{\ut_{i}}&{\cal U}^i_{\vt_{s}}&{\cal U}^i_{\ut_{s}}\\
    0&0&0&{\cal U}^s_{\vt_{i}}&{\cal U}^s_{\ut_{i}}&{\cal U}^s_{\vt_{s}}&{\cal U}^s_{\ut_{s}}\\
  \end{array}
  \right)
  \left(
  \begin{array}{c}
    \delta \pt  \\
    \deltat \vt_f \\ \deltat \ut_f \\ 
    \deltat \vt_{i} \\ \deltat \ut_{i}\\
    \deltat \vt_{s} \\ \deltat \ut_{s}
  \end{array}
  \right) = \left(
  \begin{array}{c}
    \bt_1 \\ \bt_2 \\ \bt_3 \\ \bt_4 \\ \bt_5 \\ \bt_6  \\ \bt_7
  \end{array}
  \right),
\end{equation}
where the right hand side vector $\mathbf{B}={\cal B}-{\cal
  A}(U^{(l)})$ is the Newton residual. The Jacobian shows the coupling
structure of the nonlinear problem~(\ref{system}). The indices ${\cal
  M}^f,{\cal M}^i,{\cal M}^s$ correspond to the degrees of freedom,
whether it belongs to a Lagrange node in the fluid, on the interface
or in the solid. The subnodes correspond to the dependency on the
unknown solution component, pressure, velocity and deformation, each
in the different domains. 

Three of the entries in bold letters, $\boldsymbol{{\cal
    M}^{s}_{\ut_i},{\cal M}^s_{\ut_s}}$ and $\boldsymbol{{\cal
    M}^i_{\ut_s}}$ are zero. As the deformation
gradient is expressed in terms of the velocity, see~(\ref{defgrad}),
the dependency of the solid equation on the solid's deformation does not appear.
The entry $\boldsymbol{{\cal
    M}^i_{\ut_{i}}}$ belongs to test functions $\phi$ that live on the
interface.  Thus, it contributes to both the solid equation and the
fluid equation, e.g.
\[
\langle \boldsymbol{{\cal M}^i_{\ut_{i}}}(\psi),\phi\rangle
=
\big( \frac{d}{d\ut_i} (J\sigmat_f\Ft^{-T})(\psi),\nabla\phi\big)_{\FL}
+\big( \underbrace{\frac{d}{d\ut_i}
  (\Ft\Sigmat_s)(\psi)}_{=0},\nabla\phi\big)_{\SO}, 
\]
where only the solid part will vanish, compare~(\ref{aletime:1.5}).
The remaining part belongs to the ALE map and these terms require the
highest computational effort.

Corresponding terms are found in
${\cal M}^{f}_{\ut_f},{\cal M}^{i}_{\ut_f}, {\cal M}^{f}_{\ut_i}$ and
also in ${\cal D}_{\ut_f}$ and ${\cal D}_{\ut_i}$, which are all
highlighted marked in gray. 
We will set these  matrix entries to zero and 
note once more that this is the only approximation within our
Newton-multigrid scheme. Sorting the unknowns as
$(p,\vt_f,\vt_i,\vt_s,\ut_f,\ut_i,\ut_s)$, the reduced system takes
the following  form  and reveals a block structure
\begin{equation}\label{REDUCEDMATRIX}
  \left(
  \begin{array}{cccc|c|cc}
    0            & {\cal D}_{\vt_f} & {\cal D}_{\vt_i} & 0 & \boldsymbol 
    0 &\boldsymbol  0&0\\ 
    {\cal M}^f_p & {\cal M}^f_{\vt_f} & {\cal M}^f_{\vt_i} & 0 & \boldsymbol 0 &\boldsymbol  0&0 \\
    {\cal M}^i_p & {\cal M}^i_{\vt_f} & {\cal M}^i_{\vt_i} & {\cal M}^i_{\vt_s} & \boldsymbol 0 &\boldsymbol  0&0\\
    {\cal M}^s_p & 0  & {\cal M}^s_{\vt_i} & {\cal M}^s_{\vt_s} & 0 & 0&0 \\
    \hline
    0 & 0 & 0 & 0 & {\cal E}^f_{\ut_{f}} & {\cal E}^f_{\ut_{i}}&0\\
    \hline
    0&0&{\cal U}^i_{\vt_{i}}&{\cal U}^i_{\vt_{s}}&0&{\cal U}^i_{\ut_{i}}&{\cal U}^i_{\ut_{i}}\\
    0&0&{\cal U}^s_{\vt_{i}}&{\cal U}^s_{\vt_{s}}&0&{\cal U}^s_{\ut_{i}}&{\cal U}^s_{\ut_{i}}\\
  \end{array}
  \right)
  \left(
  \begin{array}{c}
    \delta \pt  \\
    \deltat \vt_f \\
    \deltat \vt_{i} \\
    \deltat \vt_{s} \\
    \deltat \ut_f \\ 
    \deltat \ut_{i}\\
    \deltat \ut_{s}
  \end{array}
  \right) = \left(
  \begin{array}{c}
    \bt_1 \\ \bt_2 \\ \bt_3 \\ \bt_4 \\ \bt_5 \\ \bt_6  \\ \bt_7
  \end{array}
  \right).
\end{equation}
The dropped ALE derivatives (bold face zeros) are the most costly
parts in matrix 
assembly. While skipping these terms does worsen Newton convergence
rates, the overall computational time can still benefit. This has been
shown in~\cite[Section 5.2.3]{Richter2017} considering a challenging
benchmark problem with large deformation. 
This reduced linear system decomposes into three sub-steps. First, the
coupled momentum equation, living in fluid and solid domain and acting
on pressure and velocity
\begin{equation}\label{problem:1}
  \left(
  \begin{array}{cccc}
    0            & {\cal D}_{\vt_f} & {\cal D}_{\vt_i} & 0\\
    {\cal M}^f_p & {\cal M}^f_{\vt_f} & {\cal M}^f_{\vt_i} & 0\\
    {\cal M}^i_p & {\cal M}^i_{\vt_f} & {\cal M}^i_{\vt_i} & {\cal M}^i_{\vt_s}\\
    {\cal M}^s_p & 0  & {\cal M}^s_{\vt_i} & {\cal M}^s_{\vt_s}\\
  \end{array}
  \right)
  \left(
  \begin{array}{c}
    \delta \pt  \\
    \deltat \vt_f \\
    \deltat \vt_{i} \\
    \deltat \vt_{s} \\
  \end{array}
  \right) = \left(
  \begin{array}{c}
    \bt_1 \\ \bt_2 \\ \bt_3 \\ \bt_4
  \end{array}
  \right).
\end{equation}
Second, the update equation for the deformation on the interface and
within the solid domain
\begin{equation}\label{problem:2}
  \left(
  \begin{array}{cc}
    {\cal U}^i_{\ut_{i}}&{\cal U}^i_{\ut_{i}}\\
    {\cal U}^s_{\ut_{i}}&{\cal U}^s_{\ut_{i}}\\
  \end{array}
  \right)
  \left(
  \begin{array}{c}
    \deltat \ut_{i}\\
    \deltat \ut_{s}
  \end{array}
  \right) = \left(
  \begin{array}{c}
    \bt_6  \\ \bt_7
  \end{array}
  \right)-
  \begin{pmatrix}
    {\cal U}^i_{\vt_{i}}&{\cal U}^i_{\vt_{s}}\\
    {\cal U}^s_{\vt_{i}}&{\cal U}^s_{\vt_{s}}
  \end{pmatrix}
  \begin{pmatrix} \vt_i \\ \vt_s \end{pmatrix},
\end{equation}
which, as a finite element discretization of the zero-order equation
$\ut_n = \ut_{n+1}+k (1-\theta)\vt_{n-1}+k\theta\vt_n$, only
involves the mass matrix on both sides, such that this update can be
performed by one vector-addition. Finally it remains to solve for the
ALE extension equation
\begin{equation}\label{problem:3}
  {\cal E}^f_{\ut_{f}}\delta\ut_{f} = \bt_5- {\cal E}^f_{\ut_{i}}\delta\ut_{i}
\end{equation}
one simple equation, usually either a vector Laplacian or a linear
elasticity problem, see~\cite[section 5.2.5]{Richter2017}. 

The main effort lies in the momentum
equations~(\ref{problem:1}), which is still a coupled fluid-solid
problem with saddle-point character due to the incompressibility. 

Details on the derivatives appearing in~(\ref{problem:1}) are given
in~\cite{FernandezMoubachir2005,ZeeBrummelenBorst2010a,ZeeBrummelenBorst2010b}
and in~\cite[Section 5.2.2]{Richter2017} in the framework of this
work. Note however that most of these terms, including all derivatives
of the Navier-Stokes equation in direction of the fluid domain
deformation $\ut_f$ are skipped, such that the resulting fluid problem
is a weighted (due to domain deformation) variant of the Navier-Stokes
equation.

\subsection{Solution of the linear problems}

The efficient solution of the linear systems arising in Newton
approximations to nonlinear fluid-structure interaction
problems is still an open problem. Lately some progress has been done
in the direction of multigrid preconditioners for the monolithic
problem~\cite{GeeKuettlerWall2010,Richter2015,AulisaBnaBornia2018,Richter2017}. In
all these contributions it has proven to be essential to apply a
partitioning into fluid-problem and solid-problem within the smoother.
The authors of~\cite{BrummelenZeeBorst2008} analyzed a simplified
fluid-structure interaction problem and showed that a partitioned
(exact) inversion of fluid and solid problem within the multigrid
solver acts as perfect smoother with convergence rates tending to zero
on finer meshes.

We shortly present the linear algebra framework used in the software
library \emph{Gascoigne 3D}~\cite{Gascoigne3D}. We are using equal-order
finite element for all unknowns, namely pressure, velocity and
deformation such that we can block all degrees of freedom locally. The 
solution $U_h$ is written as
\[
U_h(x) = \sum_{i=1}^{N_h}\Ut_i \phi_h^{(i)}(x),\quad \Ut_i
= \begin{pmatrix}
  \pt_i \\ \vt_i \\ \ut_i
\end{pmatrix}\in \mathds{R}^{2d+1}. 
\]
By $N_h$ we denote the number of degrees of freedom (for every
unknown), by $d$ the dimension. Likewise, the system matrix $\At$ is
a matrix with block structure, i.e. $\At\in
\mathds{R}^{N_h(2d+1)\times N_h(2d+1)}$ with $\At_{ij}\in
\mathds{R}^{(2d+1)\times (2d+1)}$. Considering the approximation
scheme described in~(\ref{problem:1}),~(\ref{problem:2})
and~(\ref{problem:3}), the first problem has $n_c^{\cal M}=d+1$
components and the extension problem consists of $n_c^{\cal E}=d$
components. In general, the complete linear algebra module is acting
on general matrices and vectors  with a block structure and local
blocks of size $n_c\times n_c$ and $n_c$, respectively. 
The linear solver is designed by the following approach:
\begin{itemize}
\item[(I)] As outer iteration we employ a GMRES method. Usually very few
  ($<10$) iterations are required such that restarting strategies are
  not necessary.
\item[(II)] The GMRES solver is preconditioned by a geometric multigrid
  method in V-cycle~\cite{BeckerBraack2000a,KimmritzRichter2010}. The
  finite element mesh  of each multigrid level resolves the
  fluid-solid interface. 
\item[(III)] As smoother in the multigrid solver we use a Vanka type
  iteration which we will outline in some detail. 
\end{itemize}

The smoother for the velocity problem and the smoother for the ALE
extension problem is of Vanka type. Let ${\cal N}_h$ be the set of
degrees of freedom of the discretization on mesh level $\Omega_h$. By
${\cal P} = \{P_1,\dots,P_{n_{\cal P}}\}$ with $P_i\subset{\cal N}_h$ we
denote a partitioning of unknowns into local patches. In the most
simple case, $P_i$ includes all degrees of freedom in one element of
the mesh. Larger patches, e.g. by combining 4 adjacent elements in 2d or
8 elements in 3d are possible. 
By
$n_{\cal P}$ we denote the number of 
patches and by $n_p$ the size of each patch, which is the number of
degrees of freedom in the patch. For simplicity, we assume that all
patches in ${\cal P}$ have the same size. By ${\cal
  R}_i:\mathds{R}^N\to\mathds{R}^{n_p}$ we denote the restriction of a
global vector to the degrees of freedom in one patch, by ${\cal R}^T_i$
the prolongation. Given a block vector $\xt\in\mathds{R}^{N_hn_c}$ and
a block matrix $\At\in \mathds{R}^{N_hn_c\times N_hn_h}$ we 
denote by 
\[
\xt_i:={\cal R}_i \xt,\quad
\At_i:={\cal R}_i \At {\cal R}_i^T
\]
the restrictions to the degrees of freedom of one patch $P_i$. 
We iterate
\begin{equation}\label{vanka}
  \begin{aligned}
    d^{(l)}_h&= b_h-A_hx_h^{(l)},\\
    x_h^{(l+1)}&= x_h^{(l)} + \omega_V \sum_{P\subset\Omega_h}
    {\cal R}_i^T\At_i^{-1} {\cal R}_i d_h^{(l)},
  \end{aligned}
\end{equation}
with a damping parameter $\omega_V\approx 0.8$. This smoother can also
be considered as a domain decomposition iteration with minimal
overlap. Numerical tests have shown that this simple Jacobi coupling
is more efficient than a corresponding Gauss-Seidel iteration.

The local matrices $\At_i$ are inverted exactly using the library
\emph{Eigen}~\cite{eigenweb}. They are of  substantial
size, for $d=3$, the local matrices corresponding to the momentum
equations~(\ref{problem:1}) have dimension $108\times 108$ if  small
patches are used and $500\times 500$ if the smoother is based on the
larger patches.

\subsection{Parallelization}

Basic features of \emph{Gascoigne 3D}~\cite{Gascoigne3D} are
parallelized based on \emph{OpenMP}~\cite{openmp}. For parallelization
of the assembly of residuals and the matrix as well as application of
the Vanka smoother~(\ref{vanka}) we use a coloring of the patches
${\cal P}$ such that no collisions appear. The usual memory bottleneck
of finite element simulations will limit the parallel efficiency of
matrix vector product and Vanka smoother. We will present some data on
the parallel performance in Section~\ref{sec:parallel}.

\section{Numerical Results}\label{sec:num}

\subsection{Problem configuration}

Two different test-cases are considered to study the performance of
the discretization and the solvers that have been presented in
Sections~\ref{sec:disc} and~\ref{sec:solver}. First, we perform a
numerical study based on the 2d fsi-3 benchmark problem that has been
defined by Hron and Turek~\cite{HronTurek2006}. Second, we present a
new 3d benchmark configuration that is based on the Hron \& Turek
problem. 

\subsubsection{2d configuration}

As two dimensional configuration we solve the nonstationary 2d fsi-3
benchmark problem that has been introduced by Hron and
Turek~\cite{HronTurek2006} and since then has been revisited in
many 
contributions~\cite{HeilHazelBoyle2008,RichterWick2010} or \cite[chapter 7]{Richter2017}. We present results for this well established benchmark
problem in order to validate the discretization and to the compare the
performance of the solver with results published in literature. The
material parameters are given in Table \ref{tb:param2d} and the
parameters yield a Reynolds number (where we choose $L=\unit[0.1]{m}$
as the diameter of the cylinder) 
\[
Re_{2d} = \frac{\bar\vt\cdot L}{\nu} = 200,
\]
showing a periodic flow pattern.

\begin{table}[h]
  \begin{center}
    \begin{tabular}{l|rr}
      \toprule
      & 2d configuration & 3d configuration \\
      \midrule
      $\bar\vt$ & $\unit[2]{m\cdot s^{-1}} $                        &  $\unit[1.75]{m\cdot s^{-1}} $                   \\
      $\rho_s,\rho_f$ & $\unit[1\,000]{kg\cdot m^{-3}}$                    &  $\unit[1\,000]{kg\cdot m^{-3}}$                 \\
      $\mu_s$ & $\unit[2\cdot 10^6]{kg\cdot m^{-1}\cdot s^{-2}}$     & $\unit[2\cdot 10^6]{kg\cdot m^{-1}\cdot s^{-2}}$\\
      $\nu_f $& $\unit[0.001]{m^2\cdot s^{-1}}$                     &  $\unit[0.001]{m^2\cdot s^{-1}}$                 \\
      $\lambda_s$ & $\unit[8\cdot 10^6]{kg\cdot m^{-1}\cdot s^{-2}}$&  $\unit[8\cdot 10^6]{kg\cdot m^{-1}\cdot s^{-2}}$\\
      \bottomrule
    \end{tabular}
  \end{center}
\caption{Parameters of the benchmark problems in 2d (left) and
  3d (right).}
\label{tb:param2d}\label{tb:param3d}
\end{table}

\subsubsection{3d configuration}

\begin{figure}[t]
  \begin{center}
    \includegraphics[width=0.9\textwidth]{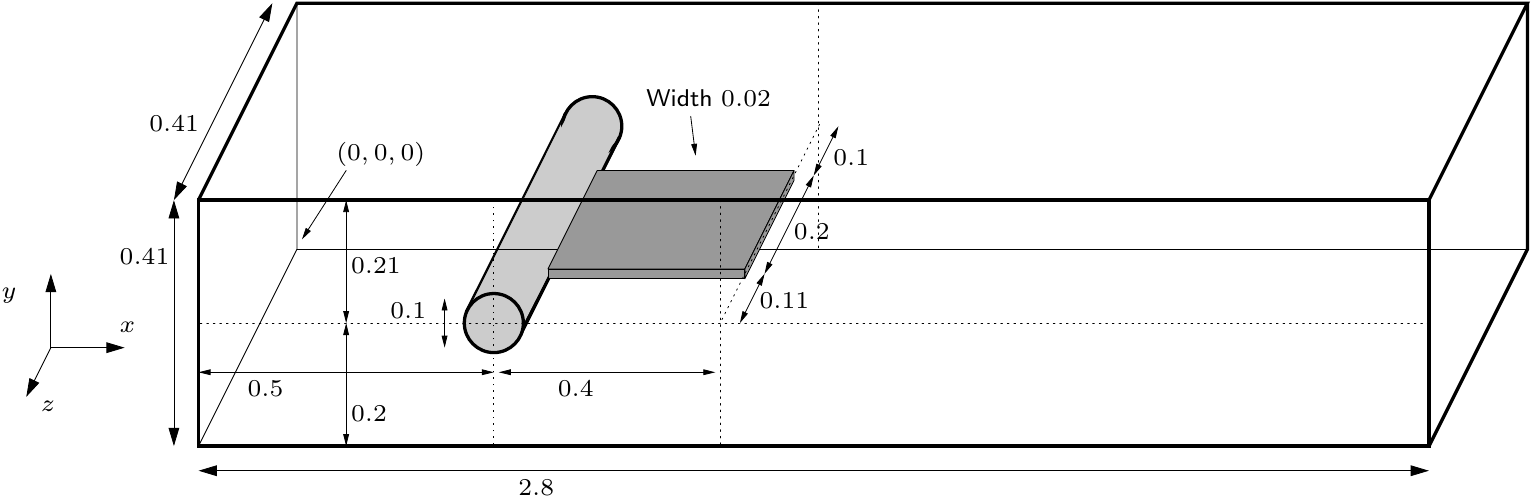}
  \end{center}
  \caption{Configuration of the 3d benchmark problem.}
  \label{fig:config3d}
\end{figure}

Figure~\ref{fig:config3d} shows the geometric configuration of the
3d benchmark problem. The computational domain with dimension $\unit[2.8]{m}\times
\unit[0.41]{m}\times \unit[0.41]{m}$ is hexahedral with a
cylinder cut out of it
\[
\begin{aligned}
  \Omega &= \{(x,y,z)\in\mathds{R}^3\,|\, 0<x<2.8,\; 0<y<0.41,\; 0<z<0.41
  \}\setminus \bar\Omega_\text{cyl},\\
  \Omega_\text{cyl} &= \{(x,y,z)\in\mathds{R}^3\,|\,
  (x-0.5)^2+(y-0.2)^2<0.05^2,\; 0<z<0.41\}.
\end{aligned}
\]
The midpoint of the cylinder is slightly non-symmetric to allow for
a stable oscillatory flow at low Reynolds numbers. Attached to the
cylinder is an elastic beam with approximate dimension $0.35\times
0.02\times 0.2$ given in initial state at time $t=0$ as 
\[
  \SO= \{(x,y,z)\in\mathds{R}^3\,|\, 0.5<x<0.9,\; 0.19<y<0.21,\;
  0.1<z<0.3\}\setminus\bar\Omega_\text{cyl} 
\]
The reference fluid domain at time $t=0$ is given by
\[  
\FL = \Omega\setminus\bar\SO.
\]

\paragraph{Boundary conditions}

The boundary of the domain is split into the \emph{inflow boundary}
$\Gamma_f^{in}$ at $x=0$, the \emph{outflow boundary}
$\Gamma_f^{out}$ at $x=2.8$, the \emph{wall boundaries} at $z=0$
and $z=0.41$ as well as $y=0$ and $y=0.41$ as well as the
\emph{cylinder boundary}  $\Gamma^{cyl}_f$ at
$(x-0.5)^2+(y-0.2)^2=0.05^2$. 
On the inflow boundary $\Gamma_f^{in}$ we prescribe a bi-parabolic
profile 
\[
\vt^{in} = \bar \vt \frac{36
  y(0.41-y)z(0.41-z)}{0.41^4}, 
\]
that satisfies $|\Gamma_f^{in}|^{-1}\int_{\Gamma_f^{in}}
\vt^{in}\,\text{d}s =\bar\vt$, 
where $\bar\vt$ is the average velocity. For regularization we suggest to
introduce a transient start-up of the inflow
\[
\vt^{in}(t) = \vt^{in} \begin{cases}
  \big(\frac{1}{2}-\frac{1}{2}\cos(\pi t)\big) &  0\le t < 1\\ 
  1 & t \ge 1.\\
\end{cases}
\]
On the remaining boundaries $\Gamma_f^{wall}\cup\Gamma_f^{cyl}$
the no-slip condition $\vt=0$ is prescribed.
For the deformation $\ut$ (both the solid deformation and the ALE
extension), a no-slip condition $\ut=0$ is prescribed 
on all boundaries. On the outer boundaries $\Gamma_\text{wall},
\Gamma_\text{in}$ and $\Gamma_\text{out}$ this condition can be
relaxed to allow for larger mesh deformations,
see~\cite[Section 5.3.5]{Richter2017}. 
\paragraph{Material Parameters}
Similar material parameters as for the 2d set are taken and the values are given in Table \ref{tb:param3d}. These parameters give a Reynolds number of
\[
Re_{3d} = \frac{\bar\vt\cdot L}{\nu} = 175,
\]
and a periodic flow pattern arises.

\subsection{Quantities of interest}

For the 2d configuration, we present the displacement at the tip of
the flag at the point $A=(0.6, 0.2)$ in x- and y-direction. In the
case of the 3d configuration we take the point $B=(0.9,0.2,0.3)$ on
the back face of the beam and 
present the displacement in x-, y- and z-direction. These values are
evaluated at every time-point. In addition we compute the drag and
lift values around the beam and cylinder. To compute the lift $\vec
f\cdot \vec e_1$ and drag forces $\vec f\cdot \vec e_2$ with $\vec e_i
= \delta_{ij}\in\mathds{R}^3$ and 
\begin{align}
  \vec f=\int_{\Gamma_f^{cyl}\cup\IN}  J\sigmat_f \Ft^{-T} n  \text{
    d}\Gamma , 
\end{align}
we evaluate the residual representation 
\begin{align*}
  \vec f_{n}&=   \big(\bar J_n (\vt_n-\vt_{n-1}),\onet_{\Gamma_{cyl}} \big)_{\FL}
    - \big( \bar J_n\bar\Ft^{-1}
    (\ut_n-\ut_{n-1})\cdot\nabla\bar\vt_n,\onet_{\Gamma_{cyl}} \big)_{\FL}\\
    &+ k A_p(U_n,\onet_{\Gamma_{cyl}}) 
    +k\theta A_F(U_n,\onet_{\Gamma_{cyl}})    
    + k(1-\theta) A_F(U_{n-1},\onet_{\Gamma_{cyl}}) \\
    &+ k(1-\theta) A_S(U_{n-1},\onet_{\Gamma_{s}}) ) 
    + k\theta A_S(U_n,\onet_{\Gamma_{s}})
\end{align*}    
where  $\onet_{\Gamma_{cyl}}$ is a finite element testfunction which is
one the cylinder $\Gamma_{cyl}$ and zero elsewhere. Thereby we can
compute the mean drag and lift value on every time interval
$I_n=[t_n,t_{n+1}]$ with very high precision. Details on the
evaluation of such surface integrals for flow problems are given
in~\cite{BraackRichter2006d} and in~\cite[Section 6.6.2]{Richter2017}
in the case of fluid-structure interactions. 



%
%




\subsection{Approximative Newton scheme (2d
  benchmark)}\label{sec:gamma}
We start by investigating the effect of the approximation of the
Jacobian in our reduced Newton scheme.  The 2d fsi-3 benchmark problem
by Hron and Turek is evaluated on the time interval $I=[5,5.5]$, where the
dynamics is fully evolved and large deformations appear. 
A similar study with the same parameters and
discretization has been performed  in \cite[chapter
  5.2.3]{Richter2017}, however, based on the full monolithic Jacobian  and
using a direct solver for the linear problems. 
The comparison
with the results in \cite{Richter2017} enables to evaluate the effects
of the presented inexact Jacobian on the Newton scheme. On the time
interval $I=[5,5.5]$ the oscillations are fully developed such that
significant oscillations appear and the geometric nonlinearities, that
come from the ALE mapping, have to be taken into account.    

\begin{figure}[t]
  \centering 
  \setlength \figureheight{3.8cm} 
  \setlength\figurewidth{0.8\textwidth}
  \small
%
%
%
\begin{tikzpicture}
\begin{axis}[%
width=0.951\figurewidth,
height=\figureheight,
at={(0\figurewidth,0\figureheight)},
scale only axis,
xmin=0,
xmax=100,
ymin=2,
ymax=20,
axis background/.style={fill=white},
title style={font=\bfseries},
title={Number of Newton iterations},
legend style={legend cell align=left,align=left,draw=white!15!black}
]
\addplot [color=mycolor1,solid,mark=o,mark options={solid}]
  table[row sep=crcr]{%
1	5\\
2	5\\
3	5\\
4	5\\
5	4\\
6	4\\
7	5\\
8	5\\
9	4\\
10	4\\
11	4\\
12	4\\
13	4\\
14	4\\
15	4\\
16	5\\
17	5\\
18	5\\
19	5\\
20	5\\
21	5\\
22	5\\
23	5\\
24	4\\
25	4\\
26	5\\
27	4\\
28	4\\
29	4\\
30	4\\
31	4\\
32	4\\
33	5\\
34	5\\
35	5\\
36	5\\
37	5\\
38	5\\
39	5\\
40	5\\
41	5\\
42	5\\
43	4\\
44	4\\
45	5\\
46	5\\
47	4\\
48	4\\
49	4\\
50	4\\
51	4\\
52	5\\
53	5\\
54	5\\
55	5\\
56	5\\
57	5\\
58	5\\
59	5\\
60	5\\
61	4\\
62	5\\
63	5\\
64	5\\
65	4\\
66	4\\
67	4\\
68	4\\
69	4\\
70	5\\
71	5\\
72	5\\
73	5\\
74	5\\
75	5\\
76	5\\
77	5\\
78	5\\
79	4\\
80	4\\
81	5\\
82	5\\
83	4\\
84	4\\
85	4\\
86	4\\
87	4\\
88	4\\
89	5\\
90	5\\
91	5\\
92	5\\
93	5\\
94	5\\
95	5\\
96	5\\
97	5\\
98	4\\
99	4\\
100	5\\
};
\addlegendentry{$\gamma\text{=0}$};

\addplot [color=mycolor2,solid,mark=square,mark options={solid}]
  table[row sep=crcr]{%
1	5\\
2	6\\
3	6\\
4	6\\
5	6\\
6	6\\
7	6\\
8	6\\
9	6\\
10	6\\
11	5\\
12	5\\
13	5\\
14	5\\
15	5\\
16	5\\
17	5\\
18	6\\
19	6\\
20	6\\
21	6\\
22	6\\
23	6\\
24	6\\
25	6\\
26	6\\
27	6\\
28	6\\
29	5\\
30	5\\
31	4\\
32	5\\
33	5\\
34	5\\
35	5\\
36	5\\
37	6\\
38	6\\
39	6\\
40	6\\
41	6\\
42	6\\
43	6\\
44	6\\
45	6\\
46	6\\
47	6\\
48	5\\
49	5\\
50	5\\
51	5\\
52	5\\
53	5\\
54	5\\
55	6\\
56	6\\
57	6\\
58	6\\
59	6\\
60	6\\
61	6\\
62	6\\
63	6\\
64	6\\
65	6\\
66	5\\
67	5\\
68	4\\
69	5\\
70	5\\
71	5\\
72	5\\
73	5\\
74	6\\
75	6\\
76	6\\
77	6\\
78	6\\
79	6\\
80	6\\
81	6\\
82	6\\
83	6\\
84	6\\
85	5\\
86	4\\
87	5\\
88	5\\
89	5\\
90	5\\
91	5\\
92	6\\
93	6\\
94	6\\
95	6\\
96	6\\
97	6\\
98	6\\
99	6\\
100	6\\
};
\addlegendentry{$\gamma\text{=0.05}$};

\addplot [color=mycolor3,solid,mark=diamond,mark options={solid}]
  table[row sep=crcr]{%
1	8\\
2	7\\
3	7\\
4	7\\
5	9\\
6	8\\
7	8\\
8	6\\
9	19\\
10	7\\
11	13\\
12	13\\
13	10\\
14	6\\
15	6\\
16	8\\
17	8\\
18	6\\
19	7\\
20	7\\
21	7\\
22	7\\
23	7\\
24	8\\
25	8\\
26	7\\
27	6\\
28	16\\
29	7\\
30	10\\
31	10\\
32	9\\
33	7\\
34	7\\
35	7\\
36	6\\
37	6\\
38	8\\
39	7\\
40	7\\
41	8\\
42	6\\
43	8\\
44	7\\
45	6\\
46	18\\
47	7\\
48	12\\
49	6\\
50	10\\
51	6\\
52	6\\
53	7\\
54	6\\
55	6\\
56	6\\
57	7\\
58	7\\
59	7\\
60	8\\
61	9\\
62	8\\
63	6\\
64	20\\
65	6\\
66	13\\
67	5\\
68	10\\
69	10\\
70	6\\
71	7\\
72	7\\
73	6\\
74	7\\
75	7\\
76	7\\
77	7\\
78	7\\
79	8\\
80	8\\
81	7\\
82	6\\
83	17\\
84	6\\
85	13\\
86	6\\
87	10\\
88	6\\
89	7\\
90	8\\
91	6\\
92	6\\
93	7\\
94	7\\
95	7\\
96	7\\
97	8\\
98	8\\
99	7\\
100	6\\
};
\addlegendentry{$\gamma\text{=0.5}$};

\end{axis}
\end{tikzpicture}%

%
%
%
\begin{tikzpicture}

\begin{axis}[%
width=0.951\figurewidth,
height=\figureheight,
at={(0\figurewidth,0\figureheight)},
scale only axis,
xmin=0,
xmax=100,
ymin=0,
ymax=8,
axis background/.style={fill=white},
title style={font=\bfseries},
title={Number of Jacobians assembled},
legend style={legend cell align=left,align=left,draw=white!15!black}
]
\addplot [color=mycolor1,solid,mark=o,mark options={solid}]
  table[row sep=crcr]{%
1	5\\
2	5\\
3	5\\
4	5\\
5	4\\
6	4\\
7	5\\
8	5\\
9	4\\
10	4\\
11	4\\
12	4\\
13	4\\
14	4\\
15	4\\
16	5\\
17	5\\
18	5\\
19	5\\
20	5\\
21	5\\
22	5\\
23	5\\
24	4\\
25	4\\
26	5\\
27	4\\
28	4\\
29	4\\
30	4\\
31	4\\
32	4\\
33	5\\
34	5\\
35	5\\
36	5\\
37	5\\
38	5\\
39	5\\
40	5\\
41	5\\
42	5\\
43	4\\
44	4\\
45	5\\
46	5\\
47	4\\
48	4\\
49	4\\
50	4\\
51	4\\
52	5\\
53	5\\
54	5\\
55	5\\
56	5\\
57	5\\
58	5\\
59	5\\
60	5\\
61	4\\
62	5\\
63	5\\
64	5\\
65	4\\
66	4\\
67	4\\
68	4\\
69	4\\
70	5\\
71	5\\
72	5\\
73	5\\
74	5\\
75	5\\
76	5\\
77	5\\
78	5\\
79	4\\
80	4\\
81	5\\
82	5\\
83	4\\
84	4\\
85	4\\
86	4\\
87	4\\
88	4\\
89	5\\
90	5\\
91	5\\
92	5\\
93	5\\
94	5\\
95	5\\
96	5\\
97	5\\
98	4\\
99	4\\
100	5\\
};
\addlegendentry{$\gamma\text{=0}$};

\addplot [color=mycolor2,solid,mark=square,mark options={solid}]
  table[row sep=crcr]{%
1	2\\
2	2\\
3	2\\
4	2\\
5	2\\
6	2\\
7	2\\
8	2\\
9	2\\
10	1\\
11	1\\
12	1\\
13	1\\
14	1\\
15	1\\
16	1\\
17	1\\
18	2\\
19	2\\
20	2\\
21	2\\
22	2\\
23	2\\
24	2\\
25	2\\
26	2\\
27	2\\
28	2\\
29	1\\
30	1\\
31	1\\
32	1\\
33	1\\
34	1\\
35	1\\
36	2\\
37	2\\
38	2\\
39	2\\
40	2\\
41	2\\
42	2\\
43	2\\
44	2\\
45	2\\
46	2\\
47	1\\
48	1\\
49	1\\
50	1\\
51	1\\
52	1\\
53	1\\
54	2\\
55	2\\
56	2\\
57	2\\
58	2\\
59	2\\
60	2\\
61	2\\
62	2\\
63	2\\
64	2\\
65	2\\
66	1\\
67	1\\
68	1\\
69	1\\
70	1\\
71	1\\
72	1\\
73	2\\
74	2\\
75	2\\
76	2\\
77	2\\
78	2\\
79	2\\
80	2\\
81	2\\
82	2\\
83	2\\
84	1\\
85	1\\
86	1\\
87	1\\
88	1\\
89	1\\
90	1\\
91	2\\
92	2\\
93	2\\
94	2\\
95	2\\
96	2\\
97	2\\
98	2\\
99	2\\
100	2\\
};
\addlegendentry{$\gamma\text{=0.05}$};

\addplot [color=mycolor3,solid,mark=diamond,mark options={solid}]
  table[row sep=crcr]{%
1	1\\
2	1\\
3	1\\
4	1\\
5	1\\
6	1\\
7	1\\
8	1\\
9	0\\
10	1\\
11	0\\
12	1\\
13	0\\
14	1\\
15	1\\
16	1\\
17	1\\
18	1\\
19	1\\
20	1\\
21	1\\
22	1\\
23	1\\
24	1\\
25	1\\
26	1\\
27	1\\
28	0\\
29	1\\
30	0\\
31	1\\
32	0\\
33	1\\
34	1\\
35	1\\
36	1\\
37	1\\
38	1\\
39	1\\
40	1\\
41	1\\
42	1\\
43	1\\
44	1\\
45	1\\
46	0\\
47	1\\
48	0\\
49	1\\
50	0\\
51	1\\
52	1\\
53	1\\
54	1\\
55	1\\
56	1\\
57	1\\
58	1\\
59	1\\
60	1\\
61	1\\
62	1\\
63	1\\
64	0\\
65	1\\
66	0\\
67	1\\
68	0\\
69	1\\
70	1\\
71	1\\
72	1\\
73	1\\
74	1\\
75	1\\
76	1\\
77	1\\
78	1\\
79	1\\
80	1\\
81	1\\
82	1\\
83	0\\
84	1\\
85	0\\
86	1\\
87	0\\
88	1\\
89	1\\
90	1\\
91	1\\
92	1\\
93	1\\
94	1\\
95	1\\
96	1\\
97	1\\
98	1\\
99	1\\
100	1\\
};
\addlegendentry{$\gamma\text{=0.5}$};

\end{axis}
\end{tikzpicture}%

%
%
%
\begin{tikzpicture}

\begin{axis}[%
width=0.951\figurewidth,
height=\figureheight,
at={(0\figurewidth,0\figureheight)},
scale only axis,
xmin=0,
xmax=100,
ymin=5,
ymax=25,
axis background/.style={fill=white},
title style={font=\bfseries},
title={Overall Computational Time},
legend style={legend cell align=left,align=left,draw=white!15!black}
]
\addplot [color=mycolor1,solid,mark=o,mark options={solid}]
  table[row sep=crcr]{%
1	20\\
2	19\\
3	19\\
4	19\\
5	15\\
6	15\\
7	19\\
8	19\\
9	15\\
10	15\\
11	16\\
12	15\\
13	15\\
14	15\\
15	15\\
16	19\\
17	19\\
18	19\\
19	19\\
20	19\\
21	19\\
22	19\\
23	19\\
24	15\\
25	15\\
26	19\\
27	15\\
28	16\\
29	16\\
30	16\\
31	16\\
32	16\\
33	19\\
34	19\\
35	19\\
36	19\\
37	19\\
38	19\\
39	19\\
40	19\\
41	19\\
42	19\\
43	15\\
44	15\\
45	19\\
46	19\\
47	15\\
48	15\\
49	15\\
50	15\\
51	15\\
52	19\\
53	19\\
54	19\\
55	19\\
56	19\\
57	19\\
58	19\\
59	19\\
60	19\\
61	15\\
62	19\\
63	19\\
64	19\\
65	16\\
66	16\\
67	16\\
68	16\\
69	16\\
70	19\\
71	19\\
72	19\\
73	19\\
74	20\\
75	19\\
76	19\\
77	19\\
78	19\\
79	15\\
80	15\\
81	19\\
82	19\\
83	15\\
84	15\\
85	15\\
86	15\\
87	15\\
88	15\\
89	19\\
90	19\\
91	19\\
92	19\\
93	19\\
94	19\\
95	19\\
96	19\\
97	19\\
98	15\\
99	15\\
100	19\\
};
\addlegendentry{$\gamma\text{=0}$};

\addplot [color=mycolor2,solid,mark=square,mark options={solid}]
  table[row sep=crcr]{%
1	11\\
2	11\\
3	11\\
4	11\\
5	11\\
6	11\\
7	11\\
8	11\\
9	11\\
10	8\\
11	7\\
12	7\\
13	7\\
14	7\\
15	7\\
16	7\\
17	7\\
18	11\\
19	11\\
20	11\\
21	11\\
22	11\\
23	11\\
24	11\\
25	11\\
26	11\\
27	11\\
28	11\\
29	7\\
30	7\\
31	6\\
32	7\\
33	7\\
34	7\\
35	7\\
36	10\\
37	11\\
38	11\\
39	11\\
40	11\\
41	11\\
42	11\\
43	11\\
44	11\\
45	11\\
46	11\\
47	7\\
48	7\\
49	7\\
50	7\\
51	7\\
52	7\\
53	7\\
54	10\\
55	11\\
56	11\\
57	11\\
58	11\\
59	11\\
60	11\\
61	11\\
62	11\\
63	11\\
64	11\\
65	11\\
66	7\\
67	7\\
68	6\\
69	7\\
70	7\\
71	7\\
72	7\\
73	10\\
74	11\\
75	11\\
76	11\\
77	11\\
78	11\\
79	11\\
80	11\\
81	11\\
82	11\\
83	11\\
84	7\\
85	7\\
86	6\\
87	7\\
88	7\\
89	7\\
90	7\\
91	10\\
92	11\\
93	11\\
94	11\\
95	11\\
96	11\\
97	11\\
98	11\\
99	11\\
100	11\\
};
\addlegendentry{$\gamma\text{=0.05}$};

\addplot [color=mycolor3,solid,mark=diamond,mark options={solid}]
  table[row sep=crcr]{%
1	10\\
2	9\\
3	9\\
4	9\\
5	16\\
6	11\\
7	9\\
8	8\\
9	14\\
10	8\\
11	9\\
12	13\\
13	7\\
14	7\\
15	9\\
16	9\\
17	14\\
18	8\\
19	9\\
20	9\\
21	9\\
22	9\\
23	9\\
24	11\\
25	11\\
26	8\\
27	8\\
28	12\\
29	8\\
30	7\\
31	11\\
32	6\\
33	10\\
34	13\\
35	13\\
36	8\\
37	9\\
38	14\\
39	9\\
40	9\\
41	11\\
42	8\\
43	9\\
44	8\\
45	8\\
46	13\\
47	8\\
48	9\\
49	7\\
50	7\\
51	8\\
52	7\\
53	8\\
54	7\\
55	8\\
56	8\\
57	9\\
58	9\\
59	9\\
60	11\\
61	16\\
62	9\\
63	8\\
64	15\\
65	8\\
66	9\\
67	7\\
68	7\\
69	11\\
70	9\\
71	8\\
72	8\\
73	8\\
74	9\\
75	9\\
76	9\\
77	9\\
78	9\\
79	11\\
80	11\\
81	8\\
82	8\\
83	12\\
84	8\\
85	9\\
86	7\\
87	7\\
88	8\\
89	13\\
90	9\\
91	7\\
92	8\\
93	14\\
94	9\\
95	9\\
96	9\\
97	11\\
98	11\\
99	8\\
100	8\\
};
\addlegendentry{$\gamma\text{=0.5}$};

\end{axis}
\end{tikzpicture}%
  \caption{Study on the effect of the non-exact Newton scheme for the
    2d benchmark problem. The Jacobian is only reassembled, if the
    Newton rate is above $\gamma$. Top: number of Newton iterations
    per time step. Middle: Number of Jacobians assembled in each time
    step. Bottom: overall computational time in each time step.}
  \label{fg:newtonsteps}
\end{figure}
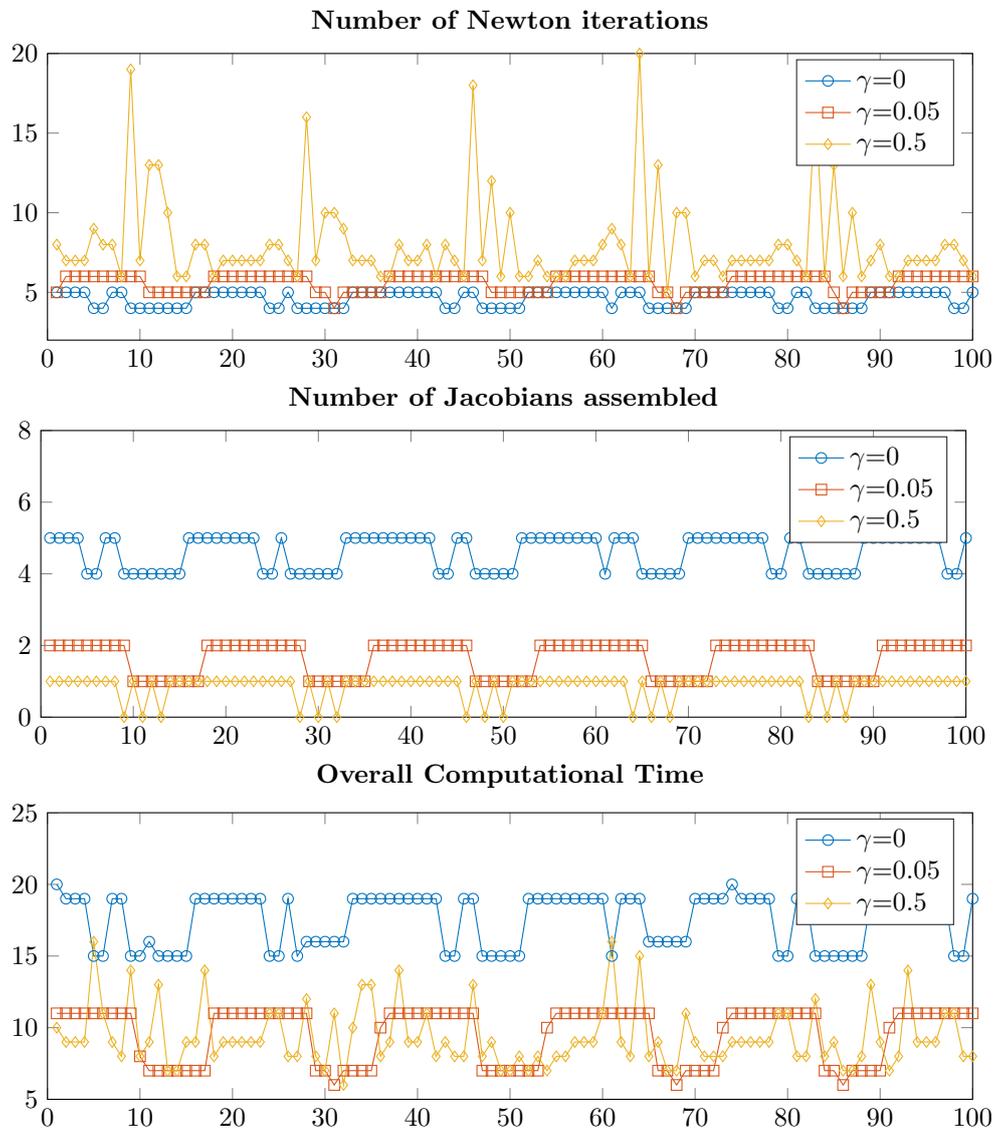 

\FloatBarrier

We only update the Jacobian of \eqref{problem:1}, the momentum
equation, if the nonlinear convergence rate, that is 
measured  as 
\begin{align}
\rho_l=\frac{\norm{{\cal B}-{\cal A}(U^{(l)})}_\infty}{\norm{{\cal B}-{\cal A}(U^{(l-1)})}_\infty},
\end{align}
is above a given threshold $\gamma_{nt}$.
The Jacobian of
\eqref{problem:3}, the mesh motion problem, is only assembled once in
the first time step, as we use a linear elasticity law. Like in
\cite{Richter2017}, we investigate the behavior for the parameters  
\begin{align*}
 \gamma_{nt}\in \{0,0.2,0.5\},
\end{align*}
where $\gamma_{nt}=0$ corresponds to the assembly of the approximated
Jacobian in every Newton step. We solve the linear systems in every
Newton step using a direct solver without any parallelization. The
computations are performed on an Intel(R) Core(TM) i7-7700 CPU @
3.60GHz. For the time stepping we use the suggested implicitly shifted
Crank-Nicolson scheme with $\theta=0.5+2k$ and the time step size
$k=0.005$s. For spatial discretization we choose equal-order
biquadratic elements  on a mesh with $80\,960$ dofs (mesh level
4). The Newton algorithm is stopped if the relative error reduces by
eight orders of magnitude ($\text{relative tol}=10^{-8}$).

In Figure \ref{fg:newtonsteps} we show the results for each time step
in the interval $I=[5,5.5]$. The top row shows that the least number
of Newton steps are required, if  $\gamma=0$ is used. This is expected as
$\gamma=0$ corresponds to the full Newton scheme that allows for
quadratic convergence. While the effect is small for $\gamma=0.05$,
the resulting Newton iteration count strongly increases for
$\gamma=0.5$, where up to 20 steps are required, compared to a limit
of 5 steps for $\gamma=0$ and 6 steps for $\gamma=0.05$. In the middle
plot of Figure~\ref{fg:newtonsteps} we give the number of Jacobians
that have to be assembled. For $\gamma=0$ these numbers obviously
correspond  to the number of Newton steps, as the Jacobian is newly
assembled in each step. For $\gamma=0.05$ and $\gamma=0.5$ the
required number of assemblies is strongly limited. Finally, the lower
plot shows the resulting computational time. Although $\gamma=0$
yields the best  convergence rates, it requires the highest
computational time. The choice $\gamma=0.05$ reduces the computational
time by a factor of 2 while still giving very robust
convergence. These results are in agreement with the study
in~\cite{Richter2017}. These results also show the large computational
time that is required for assembling the Jacobian and preparing the
multigrid smoother.

\begin{table}[b]
  \centering
  \begin{tabular}{l|rrrr}
    \toprule
    Matrix ass. tolerance & $\gamma=0.0$ & $\gamma=0.05$ & $\gamma=0.2$ & $\gamma=0.5$ \\ \midrule
    Total Newton steps & 460 & 559 & 741 &800 \\
    Jacobians assembled & 460 & 164 & 110 &85\\
    Total Time (seconds) & 1753 & 950 & 899 &936\\
    \bottomrule
  \end{tabular}
  \caption{Accumulated number of Newton steps, assemblies of the
    Jacobian in Equation \eqref{problem:1} and the total time (in
    seconds) for all 100 time steps for different values of
    $\gamma_{nt}$}
  \label{tb:accumulatedJacobians}
\end{table}

\FloatBarrier

We can see in Table \ref{tb:accumulatedJacobians}, where we collect
the accumulated numbers for the complete interval $I=[5,5.5]$ that we
need 460 Newton steps, if we assemble the Jacobian in \eqref{problem:1} in
every Newton step. As we neglect the sensitivity information with
respect to the mesh motion, we still have an inexact Newton
scheme. Nevertheless, we need less Newton steps compared to the use of and exact Jacobian as in \cite{Richter2017}, where 532
Newton steps were required for the same setting. This is in line with
the numerical tests on the inexact Jacobian for the 2d fsi-3 benchmark
results in  \cite{Richter2017}, where in first numerical studies no
disadvantages due to the inexact Jacobian could be
observed. Nevertheless, the better convergence rate is surprising. The
direct solver UMFPACK~\cite{Davis2014} has  difficulties to solve the
exact Jacobian accurately enough as reported 
in \cite{Richter2015,Richter2017}, which could be the reason for the higher number
of Newton steps. A similar study in~\cite{AulisaBnaBornia2018} shows
better robustness of the linear solver MUMPS~\cite{MUMPS:2}. 
The condition numbers for the matrices of the
subproblems \eqref{problem:1}, \eqref{problem:2} and \eqref{problem:3}
are much better then for the exact Jacobian as already analyzed in
\cite{Richter2017}. 
 
 The behavior with respect to the parameter $\gamma_{nt}$ is  comparable to the results in \cite{Richter2017}. For the pure Newton scheme a maximum of 5 Newton steps is required in comparison to 20 Newton steps for $\gamma_{nt}=0.5$. With respect to computational time, Table \ref{tb:accumulatedJacobians} shows that $\gamma=0.2$ is most efficient, as the reduced time to assemble the Jacobian and the increased time, due to more Newton steps balances best. The inexact Jacobian only has minor influence on the sensitivity of the Newton scheme with respect to the parameter $\gamma_{nt}$.

\begin{table}[h]
  \centering
  \begin{tabular}{l|cccccc}
    \toprule
    mesh level & 1 & 2 & 3 &4 &5 & 6 \\ \midrule
    dofs 2d  &1\,440	&5\,360	&20\,640	&80\,960	&320\,640	&1\,276\,155\\\hline
    dofs 3d  &63\,826	&463\,988	&3\,531\,304& -& -&-\\
    \bottomrule
  \end{tabular}
  \caption{Degrees of freedom for 2d and 3d configuration on every
    refinement level}
  \label{table:dofs}
\end{table}
 
\subsection{Reference values}
All presented solutions in the following sections are computed by
using a time stepping scheme with $k=0.004$s to compute a solution on
the time interval $I=[0,8]$ on all mesh levels indicated in
Table~\ref{table:dofs}. The corresponding solutions at time $t=8$s act
as initial values for further computations on the interval $I=[8,10]$
based on the time step sizes $k=0.004$s, $k=0.002$s and $k=0.001$s. To
avoid inaccuracies in the reference values due a rapid change of the
numerical discretization parameters, we only present results on the
interval $I=[9,10]$. A similar approach on adaptive
time-stepping schemes is demonstrated in~\cite{FailerWick} and shows
accurate results.

\FloatBarrier

\subsubsection{Reference values for the 2d configuration}

We summarized the maximal and minimal values for the functionals on
various refinement levels and time step sizes in Table
\ref{table:fsi2d_minmax_func_values}. The values indicate convergence
of the algorithm in space and a dominance of the spatial
discretization error on the coarse grids in comparison to the temporal
discretization error. These results are in very good agreement to the
values found in literature~\cite{TurekHronMadlikRazzaqWobkerAcker2010}.

\begin{table}[t]
  \centering
  \begin{tabular}{l|llll}
    \toprule
    level & \multicolumn{1}{c}{ $u_x\cdot10^{-3}$}&
    \multicolumn{1}{c}{ $u_y\cdot10^{-3}$}&\multicolumn{1}{c}{drag
      $\cdot 10^{2}$ }& \multicolumn{1}{c}{lift$\cdot 10^{2}$} \\
    \midrule
    2 & -2.5207 $\pm$ 2.4006 & 1.2285 $\pm$ 32.6701 & 4.4132 $\pm$ 0.2599 & 0.0921 $\pm$ 1.6816\\
    3 & -3.3174 $\pm$ 3.1032 & 1.2753 $\pm$ 36.8303 & 4.5564 $\pm$ 0.2941 & 0.0998 $\pm$ 1.4003\\
    4 & -2.8430 $\pm$ 2.6869 & 1.4665 $\pm$ 34.6516 & 4.5892 $\pm$ 0.2703 & 0.0363 $\pm$ 1.5581\\
    5 & -2.8716 $\pm$ 2.7174 & 1.4960 $\pm$ 34.8656 & 4.6031 $\pm$ 0.2778 & 0.0248 $\pm$ 1.5730\\
    6 & -2.8644 $\pm$ 2.7111 & 1.4995 $\pm$ 34.8329 & 4.6043 $\pm$
    0.2787 & 0.0237 $\pm$ 1.5737\\
    \multicolumn{2}{c}{}\\
    lev & \multicolumn{1}{c}{ $u_x\cdot10^{-3}$}& \multicolumn{1}{c}{
      $u_y\cdot10^{-3}$}&\multicolumn{1}{c}{drag $\cdot 10^{2}$ }&
    \multicolumn{1}{c}{lift$\cdot 10^{2}$} \\
    \midrule 
    2 & -2.6363 $\pm$ 2.5088 & 1.1688 $\pm$ 33.2886 & 4.4445 $\pm$ 0.2741 & 0.0667 $\pm$ 1.5742\\
    3 & -3.2725 $\pm$ 3.0748 & 1.2874 $\pm$ 36.7999 & 4.5753 $\pm$ 0.2964 & 0.0683 $\pm$ 1.3963\\
    4 & -2.8466 $\pm$ 2.6874 & 1.4604 $\pm$ 34.6813 & 4.5915 $\pm$ 0.2702 & 0.0319 $\pm$ 1.5509\\
    5 & -2.8850 $\pm$ 2.7255 & 1.4774 $\pm$ 34.9795 & 4.6037 $\pm$ 0.2786 & 0.0252 $\pm$ 1.5675\\
    6 & -2.8841 $\pm$ 2.7250 & 1.4785 $\pm$ 34.9845 & 4.6050 $\pm$
    0.2798 & 0.0242 $\pm$ 1.5699\\
    \multicolumn{2}{c}{}\\
    lev & \multicolumn{1}{c}{ $u_x\cdot10^{-3}$}& \multicolumn{1}{c}{
      $u_y\cdot10^{-3}$}&\multicolumn{1}{c}{drag $\cdot 10^{2}$ }&
    \multicolumn{1}{c}{lift$\cdot 10^{2}$} \\
    \midrule
    2 & -2.7866 $\pm$ 2.6462 & 1.1851 $\pm$ 33.9983 & 4.4712 $\pm$ 0.2887 & 0.0439 $\pm$ 1.4837\\
    3 & -3.2432 $\pm$ 3.0478 & 1.2869 $\pm$ 36.7179 & 4.5884 $\pm$ 0.2979 & 0.0531 $\pm$ 1.4114\\
    4 & -2.8317 $\pm$ 2.6716 & 1.4550 $\pm$ 34.6089 & 4.5925 $\pm$ 0.2686 & 0.0297 $\pm$ 1.5425\\
    5 & -2.8844 $\pm$ 2.7234 & 1.4674 $\pm$ 34.9896 & 4.6034 $\pm$ 0.2775 & 0.0250 $\pm$ 1.5610\\
    6 & -2.8900 $\pm$ 2.7290 & 1.4690 $\pm$ 35.0322 & 4.6049 $\pm$
    0.2791 & 0.0245 $\pm$ 1.5659\\
    \bottomrule
  \end{tabular}
  \caption{Results of the 2d fsi-3 Benchmark with time step size
    $k=0.004$s, $k=0.002s$ and $k=0.001$s.}
  \label{table:fsi2d_minmax_func_values}
\end{table}

\subsubsection{Reference values for the 3d configuration}

\begin{figure}[h]
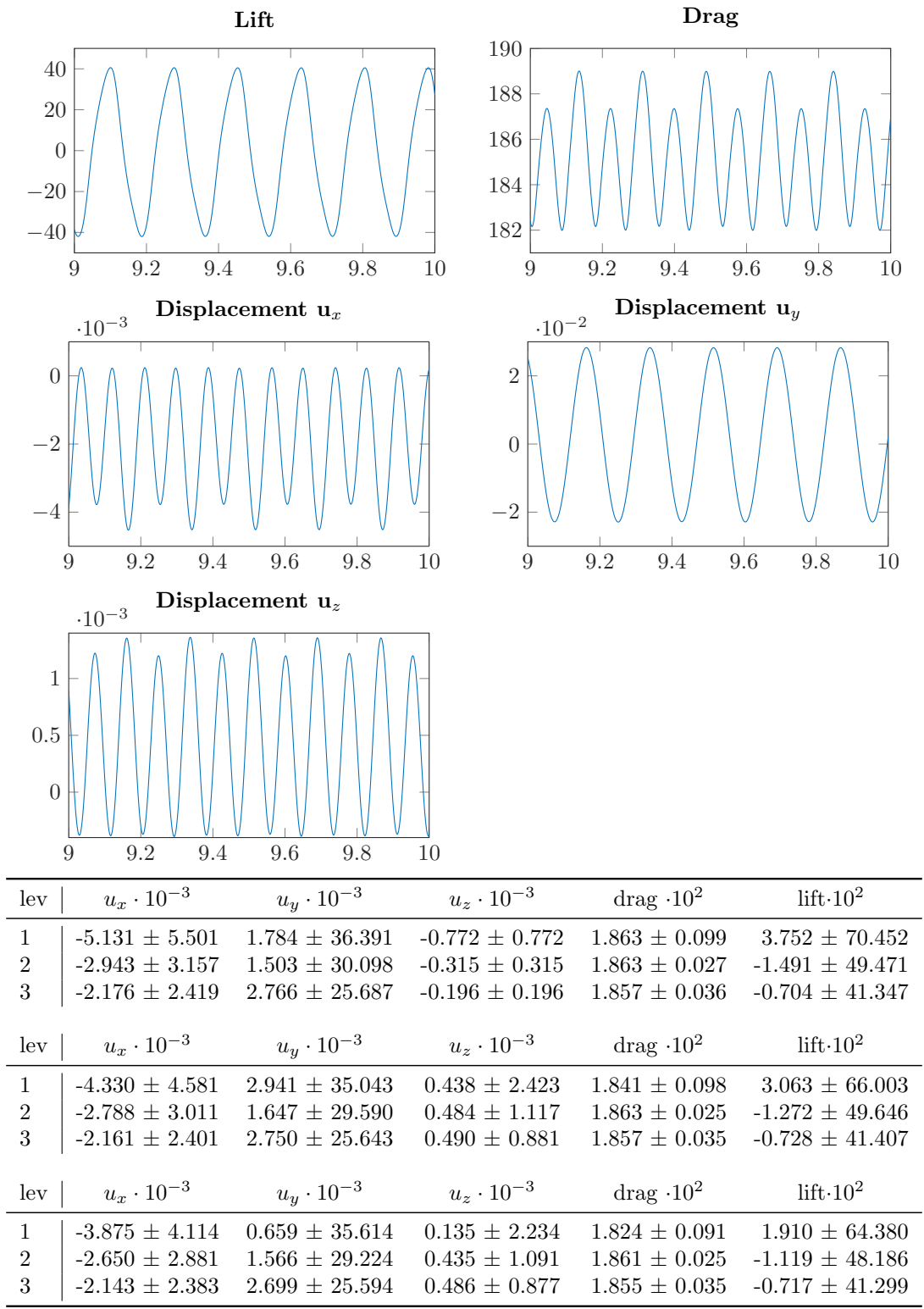

\centering
  \setlength \figureheight{0.2\textwidth} 
  \setlength\figurewidth{0.35\textwidth} 
  \small
  \begin{tabular}{cc}
    \input{fsi3d_lift.tikz}&
    \input{fsi3d_drag.tikz}\\
    \input{fsi3d_ux.tikz}&
    \input{fsi3d_uy.tikz}\\
    \input{fsi3d_uz.tikz}&
  \end{tabular}

  \begin{tabular}{l|ccccc}
    \toprule
    lev & \multicolumn{1}{c}{ $u_x\cdot10^{-3}$}& \multicolumn{1}{c}{
      $u_y\cdot10^{-3}$}& \multicolumn{1}{c}{
      $u_z\cdot10^{-3}$}&\multicolumn{1}{c}{drag $\cdot 10^{2}$ }&
    \multicolumn{1}{c}{lift$\cdot 10^{2}$} \\
    \midrule
    1 & -5.131 $\pm$ 5.501 & 1.784 $\pm$ 36.391 & -0.772 $\pm$ 0.772 & 1.863 $\pm$ 0.099 & \phantom{-}3.752 $\pm$ 70.452\\
    2& -2.943 $\pm$ 3.157 & 1.503 $\pm$ 30.098 & -0.315 $\pm$ 0.315 & 1.863 $\pm$ 0.027 & -1.491 $\pm$ 49.471\\
    3 & -2.176 $\pm$ 2.419 & 2.766 $\pm$ 25.687 & -0.196 $\pm$ 0.196 &
    1.857 $\pm$ 0.036 & -0.704 $\pm$ 41.347\\
    \multicolumn{4}{c}{}\\
    lev & \multicolumn{1}{c}{ $u_x\cdot10^{-3}$}& \multicolumn{1}{c}{ $u_y\cdot10^{-3}$}& \multicolumn{1}{c}{ $u_z\cdot10^{-3}$}&\multicolumn{1}{c}{drag $\cdot 10^{2}$ }& \multicolumn{1}{c}{lift$\cdot 10^{2}$} \\\midrule
    1 & -4.330 $\pm$ 4.581 & 2.941 $\pm$ 35.043 & 0.438 $\pm$ 2.423 & 1.841 $\pm$ 0.098 & \phantom{-}3.063 $\pm$ 66.003\\
    2 & -2.788 $\pm$ 3.011 & 1.647 $\pm$ 29.590 & 0.484 $\pm$ 1.117 & 1.863 $\pm$ 0.025 & -1.272 $\pm$ 49.646\\
    3 & -2.161 $\pm$ 2.401 & 2.750 $\pm$ 25.643 & 0.490 $\pm$ 0.881 &
    1.857 $\pm$ 0.035 & -0.728 $\pm$ 41.407\\
    \multicolumn{4}{c}{}\\
    lev & \multicolumn{1}{c}{ $u_x\cdot10^{-3}$}& \multicolumn{1}{c}{
      $u_y\cdot10^{-3}$}& \multicolumn{1}{c}{
      $u_z\cdot10^{-3}$}&\multicolumn{1}{c}{drag $\cdot 10^{2}$ }&
    \multicolumn{1}{c}{lift$\cdot 10^{2}$} \\
  \midrule
  1 & -3.875 $\pm$ 4.114 & 0.659 $\pm$ 35.614 & 0.135 $\pm$ 2.234 & 1.824 $\pm$ 0.091 & \phantom{-}1.910 $\pm$ 64.380\\
  2 & -2.650 $\pm$ 2.881 & 1.566 $\pm$ 29.224 & 0.435 $\pm$ 1.091 & 1.861 $\pm$ 0.025 & -1.119 $\pm$ 48.186\\
  3 & -2.143 $\pm$ 2.383 & 2.699 $\pm$ 25.594 & 0.486 $\pm$ 0.877 &
  1.855 $\pm$ 0.035 & -0.717 $\pm$ 41.299\\
  \bottomrule
  \end{tabular}
  
  \caption{3d fsi-3 configuration. Top:  functional values as function
    over the time interval $I=[9,10]$ for $k=0.001$s and mesh level
    3. Bottom: results for time step sizes  $k=0.004$s, $k=0.002$s,
    $k=0.001$s and all three mesh levels.}
  \label{fig:fsi3d_func_values}
\end{figure}

In 3d, we evaluate the displacement of the elastic beam in the point
$B$ and also compute  the drag and lift coefficients around the whole
cylinder and the flag. In Figure \ref{fig:fsi3d_func_values} we show
the different 
functionals as function over the time 
interval $I=[9,10]$. In addition, we summarized the maximal and
minimal value for different meshes and for different time step
sizes $k$. To draw a conclusion on the convergence or to present
reference values, the computation has to be repeated on even finer meshes
in the future. 

\FloatBarrier

\subsection{Performance of the linear solver}

To test the linear iterative solver presented in Section
\ref{sec:solver}, we recomputed the solution on different mesh levels
for the 2d and 3d benchmark configuration on the time interval
$I=[9,9.5]$ with time-step size $k=0.002$s (250 steps). The beam
oscillates in this time interval. Hence, due to the strong coupling,
the solution of the Newton system is very challenging and the fluid as
well as the solid elasticity problem have both to be solved very accurately.

The Newton algorithm in every time step terminates, if the residual is
reduced by eight orders of magnitude ($\text{relative tol}=10^{-8}$)
or if the absolute value, so the residual, falls below $10^{-8}$. In
every Newton step, the iterative solver for the linear problem
\eqref{problem:1} reduces the error by a factor of $10^{-4}$. The
parameter $\gamma=0.05$ is chosen as in Section \ref{sec:gamma} to
decide, if the Jacobian of the momentum equation \eqref{problem:1} is
reassembled in the next Newton step. The mesh motion subproblem
\eqref{problem:3} is a linear elasticity problem and hence can be
solved very efficiently with the geometric multigrid
solver. Nevertheless, as we have to solve it after every Newton step, the solution of the linear system has still a high contribution to the computational time. The matrix for the linear meshmotion problem \eqref{problem:3} only has to be assembled once in the first step. 

In the following, we will only present averaged values. By ``mean time
per Newton step'' we denote the average time of each step, measured
over all 250 time steps.  Hence, this average value also includes the
time to reassemble the Jacobian, whose assembly incidence depends on
the Newton rate, see Section~\ref{sec:gamma}. To make the values
comparable with other solution approaches, we additionally present  the
mean time to assemble one Jacobian of the momentum equation \eqref{problem:1}. In the case
of the direct solver, this includes the times for preparation and
computation of the LU decomposition. In the case of the ILU and Vanka
smoother the assemble times include the time to compute the ILU or the
LU of the block matrices $\At_i$.

\subsubsection{Dependency on the Vanka patch size (2d fsi-3)}

Concerning the Vanka smoother, the question arises, how large we should
choose the patches $P_i$ to solve the linear system coming from the
momentum equation
\eqref{problem:1} most efficiently. The simple structure of the Vanka
solver enables to use different patch sizes in the fluid and solid
domain. To test different blocking strategies we recorded the
computational time for the 2d fsi-3 benchmark on the finest mesh level 6
and present the mean number of Newton steps and matrix assemblies per
time step in Table \ref{table:testVankaBlock}. We either choose
patches consisting of one element ($n_p=3^2\cdot 3=27$) or patches
stretching over four adjacent elements ($n_p=5^2\cdot 3=75$). This
yields local matrices of size $\At_i\in\mathds{R}^{27\times 27}$ or
$\At_i\in\mathds{R}^{75\times 75}$ if larger patches are used.

We can observe that the minimal number of GMRES steps to solve
\eqref{problem:1} in every Newton step can be obtained, by using
$n_p=75$. If we only use the degrees of freedom of one element as
block on the solid domain, the number of GMRES steps increases and the
Newton convergence suffers. This effect cannot be observed, if we only
use smaller patches within the fluid domain, but large patches in the
solid. As computational times are reasonable small for 2d
computations, we will always use larger patches of size $n_p=75$ in 
the Vanka smoother.

In 3d the same blocking strategy would correspond to combining 8 elements to
one block, resulting in $n_p = 5^3\cdot 4 = 500$ and matrices of size
$\At_i\in\mathds{R}^{500\times 500}$. This strategy is forbiddingly
expensive with increasing  memory and time consumption for each
block-LU. As the results in Table \ref{table:testVankaBlock} show that
it is sufficient to use small patches in the fluid domain, we will
combine large patches with $n_p=500$ in the solid with smaller patches
of size $n_p = 3^3\cdot 4=108$ in the fluid domain for all 3d
computations to follow. 

\begin{table}[t]
\centering
\begin{tabular}{l|rrrr}
  \toprule
 $n_P$	& $\FL$: 27/$\SO$: 27 & $\FL$: 75/$\SO$: 75&$\FL$: 75/$\SO$:
  27& $\FL$: 27/$\SO$: 75\\
  \midrule
 Newton steps     	&6.87& 5.10&6.83 &5.10\\
 Matrix assemblies	&2.87&1.23&2.86&1.23\\
 GMRES per Newton	&20.58&12.60&17.16&15.26\\
 Relative comp. time &$100\%$&$55\%$&$95\%$&$61\%$\\
 \bottomrule
\end{tabular}
\caption{Vanka Blocking strategy for the 2d test case. We either
  choose every element as 
  one block or combine 4 elements to one block on the fluid ($\FL$)
  and solid ($\SO$) domain. We present the number of Newton steps,
  matrix assembles and number of GMRES steps per linear solve of
  \eqref{problem:1} on mesh level 6 and the computational time
  relative to the $\FL$: 1/$\SO$: 1 case}
\label{table:testVankaBlock} 
\end{table}

\subsubsection{Geometric multigrid performance in 2d and
  3d (sequential computations)}\label{sec:num-geo-mult} 

\begin{figure}[t]
  \centering
  \setlength\figureheight{0.28\textwidth} 
  \setlength\figurewidth{0.45\textwidth} 
  \small

  \begin{tabular}{cc}
    \multicolumn{2}{c}{\textbf{Average memory usage}}\\
    \begin{tikzpicture}
\begin{loglogaxis}[%
width=\figurewidth,
height=\figureheight,
scale only axis,
x label style={anchor=north, below=-12mm},
xlabel=dofs,
y label style={below=-3mm},
ylabel=kB,
title style={font=\bfseries},
legend style={at={(0,\figureheight)},xshift=0.2cm,anchor=north west,nodes=right}] 
]
\addplot [
color=mycolor1,
solid,
mark=*
]
table[row sep=crcr]{
5360	2.49E+04
\\
20640	1.24E+05
\\
80960	5.54E+05
\\
320640	2.55E+06
\\
};
\addlegendentry{direct}
\addplot [
color=mycolor2,
solid,
mark=square*
]
table[row sep=crcr]{
5360	1.16E+04
\\
20640	3.68E+04
\\
80960	1.37E+05
\\
320640	5.36E+05
\\
1276155 2.10E+06
\\
};
\addlegendentry{ilu}
\addplot [
color=mycolor3,
solid,
mark=diamond*
]
table[row sep=crcr]{
5360	1.47E+04
\\
20640	5.14E+04
\\
80960	1.97E+05
\\
320640	7.81E+05
\\
1276155 3.12E+06
\\
};
\addplot [
color=mycolor5,
dashed,
]
table[row sep=crcr]{
5360	0.5E+04\\
3531304	120E+04\\
};
\addlegendentry{lin}
\addlegendentry{vanka}
\end{loglogaxis}
\end{tikzpicture}%
&
    \begin{tikzpicture}
\begin{loglogaxis}[%
width=\figurewidth,
height=\figureheight,
scale only axis,
x label style={anchor=north, below=-12mm},
xlabel=dofs,
ymin=5e5,
ymax=1e8,
title style={font=\bfseries},
legend style={at={(0,\figureheight)},xshift=0.2cm,anchor=north west,nodes=right}] 
]
\addplot [
color=mycolor1,
solid,
mark=*
]
table[row sep=crcr]{
63826	1.62E+06\\
463988	3.40E+07\\
};
\addlegendentry{direct}
\addplot [
color=mycolor2,
solid,
mark=square*
]
table[row sep=crcr]{
63826	6.22E+05\\
463988	2.95E+06\\
3531304	2.14E+07\\
};
\addlegendentry{ilu}
\addplot [
color=mycolor3,
solid,
mark=diamond*
]
table[row sep=crcr]{
63826	7.03E+05\\
463988	3.59E+06\\
3531304	2.65E+07\\
};
\addlegendentry{vanka}
\addplot [
color=mycolor5,
dashed,
]
table[row sep=crcr]{
63826	10E+05\\
3531304	500E+05\\
};
\addlegendentry{lin}
\end{loglogaxis}
\end{tikzpicture}%
 \\[2mm]
    \multicolumn{2}{c}{\textbf{Mean time per Newton step}}\\\\
    \begin{tikzpicture}
\begin{loglogaxis}[%
width=\figurewidth,
height=\figureheight,
scale only axis,
x label style={anchor=north, below=-12mm},
xlabel=dofs,
y label style={below=-3mm},
ylabel=seconds,
title style={font=\bfseries},
legend style={at={(0,\figureheight)},xshift=0.2cm,anchor=north west,nodes=right}] 
]
\addplot [
color=mycolor1,
mark=*
]
table[row sep=crcr]{
5360	0.059202501954652\\
20640	0.263279445727483\\
80960	1.45988363072149\\
320640	9.52222563315426\\
};
\addlegendentry{direct}
\addplot [
color=mycolor2,
mark=square*
]
table[row sep=crcr]{
5360	0.07062353858145\\
20640	0.326705426356589\\
80960	1.53010886469673\\
320640	6.72896907216495\\
1276155 32.2277975270479\\
};
\addlegendentry{ilu}
\addplot [
color=mycolor3,
mark=diamond*
]
table[row sep=crcr]{
5360	0.078347689898199\\
20640	0.418324282389449\\
80960	2.01619751166407\\
320640	9.06886996904025\\
1276155 43.8754858934169\\
};
\addlegendentry{vanka}
\addplot [
color=mycolor5,
dashed
]
table[row sep=crcr]{
5360	0.025\\
1276155 6.5\\
};
\addlegendentry{lin}
\end{loglogaxis}
\end{tikzpicture}%
&
    \begin{tikzpicture}
\begin{loglogaxis}[%
width=\figurewidth,
height=\figureheight,
scale only axis,
x label style={anchor=north, below=-12mm},
xlabel=dofs,
title style={font=\bfseries},
legend style={at={(0,\figureheight)},xshift=0.2cm,anchor=north west,nodes=right}] 
]
\addplot [
color=mycolor1,
mark=*
]
table[row sep=crcr]{
63826	11.6106595905989\\
463988	1121.63969673995\\
};
\addlegendentry{direct}
\addplot [
color=mycolor2,
mark=square*
]
table[row sep=crcr]{
63826	5.04457840236686
\\
463988	34.2643584229391
\\
3531304	401.408582317073
\\
};
\addlegendentry{ilu}
\addplot [
color=mycolor3,
mark=diamond*
]
table[row sep=crcr]{
63826	4.16498852333589
\\
463988	36.4683651642475
\\
3531304	435.690038819876
\\
};
\addlegendentry{vanka}
\addplot [
color=mycolor5,
dashed
]
table[row sep=crcr]{
63826	2
\\
463988	16
\\
3531304	128
\\
};
\addlegendentry{lin}
\end{loglogaxis}
\end{tikzpicture}%
\\[2mm]
    \multicolumn{2}{c}{\textbf{Time for one matrix assembly}}\\
    \begin{tikzpicture}
\begin{loglogaxis}[%
width=\figurewidth,
height=\figureheight,
scale only axis,
x label style={anchor=north, below=-12mm},
xlabel=dofs,
y label style={below=-3mm},
ylabel=seconds,
title style={font=\bfseries},
legend style={at={(0,\figureheight)},xshift=0.2cm,anchor=north west,nodes=right}] 
]
\addplot [
color=mycolor1,
mark=*
]
table[row sep=crcr]{
5360	0.094405286343612\\
20640	0.469456066945607\\
80960	3.53991228070175\\
320640	26.4057192982456\\
};
\addlegendentry{direct}
\addplot [
color=mycolor2,
mark=square*
]
table[row sep=crcr]{
5360	0.031531531531532\\
20640	0.125144032921811\\
80960	0.52590308370044\\
320640	2.18950819672131\\
1276155 9.21297872340425\\
};
\addlegendentry{ilu}
\addplot [
color=mycolor3,
mark=diamond*
]
table[row sep=crcr]{
5360	0.029561403508772\\
20640	0.110334728033473\\
80960	0.442202643171806\\
320640	1.754250000000000\\
1276155 7.01077669902913\\
};
\addlegendentry{vanka}
\addplot [
color=mycolor5,
dashed,
]
table[row sep=crcr]{
5350	0.015\\
3531304	3.6\\
};
\addlegendentry{lin}
\end{loglogaxis}
\end{tikzpicture}%
 &
    \begin{tikzpicture}
\begin{loglogaxis}[%
width=\figurewidth,
height=\figureheight,
scale only axis,
x label style={anchor=north, below=-12mm},
xlabel=dofs,
title style={font=\bfseries},
legend style={at={(0,\figureheight)},xshift=0.2cm,anchor=north west,nodes=right}] 
]
\addplot [
color=mycolor1,
mark=*
]
table[row sep=crcr]{
63826	43.3344841269841\\
463988	5487.9881092437\\
};
\addlegendentry{direct}
\addplot [
color=mycolor2,
mark=square*
]
table[row sep=crcr]{
63826	5.89983766233766\\
463988	22.4867213114754\\
3531304	158.961139705882\\
};
\addlegendentry{ilu}
\addplot [
color=mycolor3,
mark=diamond*
]
table[row sep=crcr]{
63826	6.33838709677419\\
463988	20.5842259414226\\
3531304	135.424741035857\\
};
\addlegendentry{vanka}
\addplot [
color=mycolor5,
dashed,
]
table[row sep=crcr]{
63826	10\\
3531304	500\\
};
\addlegendentry{lin}
\end{loglogaxis}
\end{tikzpicture}%
  \end{tabular}
  \caption{Performance of the multigrid solver in 2d (left) and 3d
    (right) in comparison to a direct solver. On different meshes with
    increasing numbers of degrees of freedom, we compare the
    performance of the direct solver UMFPACK~\cite{Davis2014} with the
    multigrid solver based on ILU smoothing and based on Vanka
    smoothing. No parallelization is employed. From top to
    bottom: average memory usage, average time per Newton step and
    time to assemble one system matrix including preparation of the
    direct solver and the smoothers. Note that we do not reassemble
    the Jacobian in every Newton step, and therefore, assembly times
    can be higher than mean Newton times (which include the
    assembly). }
  \label{fig:memoryconsumtion}
\end{figure}
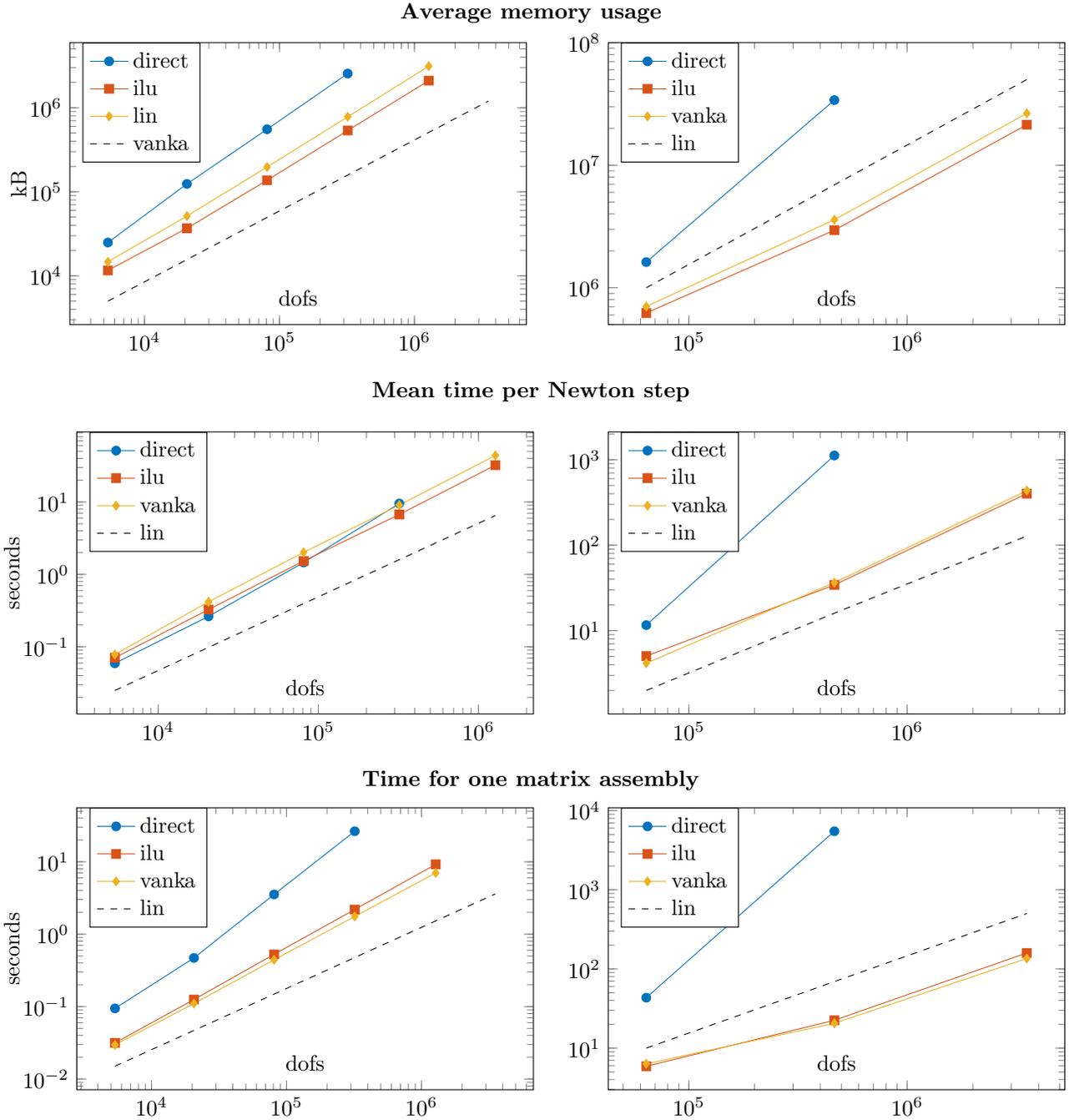

All computations have been carried out on an Intel(R) Core(TM) i7-7700
CPU @ 3.60GHz. Single Core performance only is used in this section. 
In Figure \ref{fig:memoryconsumtion} we show the results for both 2d
and 3d benchmark problems on sequences of meshes.

In the top row we present the memory consumption (in 3d, the finest
mesh level exceeded the available memory). In particular the 3d
results show the expected
superiority of iterative solvers as compared to the 
direct linear Solver UMFPACK~\cite{Davis2014} with a non-optimal
scaling. The Vanka smoother requires slightly more memory which comes
from the overlap of degrees of freedom between the different blocks.
The middle plot of Figure~\ref{fig:memoryconsumtion} shows the
resulting computational time. According to our previous
study~\cite{Richter2015}, the multigrid method is not able to beat the
direct solver in 2d. The situation dramatically changes in 3d, where
the direct solver shows a strongly non-optimal scaling. The multigrid
solvers shows nearly linear scaling for both ILU and Vanka
smoothing. Concerning the ILU smoother, this is an improvement to our
previous study presented in~\cite{Richter2015}, where the multigrid
solver was performed in a purely monolithic setting and an ILU that
consists of local blocks coupling pressure, velocity and
deformation. Here no convergence could be achieved on fine meshes. We
note that the 3d benchmark problem considered in this paper is by far
more challenging than the problem investigated
in~\cite{Richter2015,Richter2017} as it comprises very large
deformation and hence strong nonlinearities in the solid and also in
the ALE map.
The lowest row shows the time for one assembly of the
Jacobian, including the computational times for preparing  the direct
solver, the ILU smoother and the Vanka smoother. Here, the main
discrepancies between the direct solver and the multigrid methods
arise. Since we do not recompute the Jacobian in every Newton step
(not even in every time step), it is no inconsistency that the assembly
time is larger than the complete time per Newton step. In 2d the
results appear slightly sub-optimal. This is due to the necessity to
assemble the matrices along the complete multigrid hierarchy yielding
a scaling of order ${\cal O}(n\log n)$.

\FloatBarrier

\begin{table}[h]
  \centering
  \begin{tabular}{l|c@{ }c@{ }c|c@{ }c@{ }c|c@{ }c@{ }c}
    \toprule
    & \multicolumn{3}{c}{ mesh level 4}& \multicolumn{3}{c}{ mesh
      level 5}& \multicolumn{3}{c}{ mesh level 6}\\
    & direct & ILU & Vanka& direct & ILU & Vanka& direct & ILU & Vanka\\
    \midrule
    Newton steps     	&5.15 &5.04 &5.14	&5.21 &5.04 &5.16
    &- &5.17 &5.10\\
    Matrix assemblies	&1.14 &0.90 &0.90	&1.14 &0.97 &0.96
    &- &0.94 &1.23\\
    GMRES per Newton	&- &11.07 &9.53	&- &11.07 &10.65	&-
    &13.08 &12.60\\
    \bottomrule
  \end{tabular}

  \bigskip
  \begin{tabular}{l|c@{ }c@{ }c|c@{ }c@{ }c|c@{ }c@{ }c}
    \toprule
    & \multicolumn{3}{c}{mesh level 1}& \multicolumn{3}{c}{mesh level
      2}& \multicolumn{3}{c}{mesh level 3}\\
    &direct&ILU&Vanka&direct&ILU&Vanka&direct&ILU&Vanka\\
    \midrule
    Newton steps     	&5.27&5.40&5.22&5.27&5.58&5.23&-&5.24&5.15\\
    Matrix assemble	&1.00&1.23&0.99&0.95&1.22&0.95&-&1.22&1.00\\
    GMRES per Newton	&-&13.20&4.81&-&14.52&9.35&-&15.33&10.59\\
    \bottomrule
  \end{tabular}

  \caption{Average number of Newton steps, matrix assemblies per time
    step and average number of GMRES steps within each Newton
    step. Top: 2d benchmark. Bottom: 3d benchmark.}
  \label{table:Newtonsteps2d3d}
\end{table}

\begin{figure}[h]
  \centering
  \setlength \figureheight{4cm} 
  \setlength\figurewidth{0.8\textwidth}
  \small
%
%
%
%
\begin{tikzpicture}

\begin{axis}[%
width=\figurewidth,
height=\figureheight,
scale only axis,
separate axis lines,
every outer x axis line/.append style={darkgray!60!black},
every x tick label/.append style={font=\color{darkgray!60!black}},
xmin=0,
xmax=250,
every outer y axis line/.append style={darkgray!60!black},
every y tick label/.append style={font=\color{darkgray!60!black}},
ymin=9,
ymax=18,
title style={font=\bfseries},
title={\text{GMRES steps per Newton step to solve the momentum equation}},
legend style={draw=darkgray!60!black,fill=white,legend cell align=left}
]
\addplot [
color=mycolor1,
mark=o,
solid
]
table[row sep=crcr]{
1 16.6666666666667\\
2 16.4\\
3 15.6\\
4 15.4\\
5 15.6\\
6 15\\
7 15.1666666666667\\
8 15.1666666666667\\
9 15.1666666666667\\
10 14.3333333333333\\
11 15.1428571428571\\
12 15\\
13 14.4\\
14 14.2\\
15 14.2\\
16 14.2\\
17 14\\
18 14.8\\
19 14\\
20 14.2\\
21 14.2\\
22 14.4\\
23 14.6\\
24 14.8\\
25 14.8\\
26 14.8\\
27 15\\
28 15.4\\
29 15.2\\
30 15.4\\
31 15.8\\
32 16\\
33 16\\
34 16.2\\
35 16.4\\
36 15.2\\
37 15.6\\
38 16.1666666666667\\
39 16.6666666666667\\
40 16.8\\
41 16.6\\
42 16.8\\
43 17.4\\
44 17.8\\
45 17.3333333333333\\
46 16.4\\
47 15.8\\
48 15.8\\
49 15.4\\
50 14.6\\
51 15\\
52 15.1666666666667\\
53 15\\
54 14.6666666666667\\
55 14.3333333333333\\
56 15\\
57 14.3333333333333\\
58 14.2\\
59 14.2\\
60 14.6\\
61 14.4\\
62 14.4\\
63 14\\
64 14.2\\
65 14.2\\
66 14.6\\
67 14.8\\
68 15\\
69 15.2\\
70 15.2\\
71 15\\
72 15.4\\
73 15.2\\
74 15.2\\
75 15.8\\
76 15.8\\
77 16\\
78 16.2\\
79 15.5\\
80 15.6666666666667\\
81 15.6\\
82 16\\
83 16\\
84 17.4\\
85 16.8\\
86 17.2\\
87 17.4\\
88 17.6\\
89 17.5\\
90 16.8\\
91 15.8\\
92 15.8\\
93 15.8\\
94 15\\
95 15.1666666666667\\
96 15\\
97 15.3333333333333\\
98 14.6666666666667\\
99 14.5\\
100 15.1428571428571\\
101 14.4\\
102 14.2\\
103 14.2\\
104 14.4\\
105 14.2\\
106 14\\
107 13.8\\
108 14.2\\
109 14.2\\
110 14.4\\
111 14.4\\
112 14.6\\
113 14.8\\
114 14.4\\
115 14.8\\
116 15.2\\
117 15\\
118 15.2\\
119 15.4\\
120 16\\
121 16\\
122 16.2\\
123 16.4\\
124 15.6666666666667\\
125 15.4\\
126 16\\
127 15.6\\
128 17\\
129 16.4\\
130 17\\
131 17\\
132 17.8\\
133 17.8\\
134 16.8\\
135 15.8\\
136 15.8\\
137 15.8\\
138 15\\
139 14.8\\
140 15.3333333333333\\
141 15\\
142 14.6666666666667\\
143 14.3333333333333\\
144 15.1428571428571\\
145 14.1666666666667\\
146 14.5\\
147 14.4\\
148 14.4\\
149 14.4\\
150 14.6\\
151 14.4\\
152 14\\
153 14.2\\
154 14.6\\
155 14.8\\
156 14.8\\
157 15\\
158 15.2\\
159 15.2\\
160 15.2\\
161 15\\
162 15.2\\
163 15.4\\
164 15.8\\
165 16\\
166 15.8\\
167 16.4\\
168 15.6666666666667\\
169 15.4\\
170 16\\
171 15.6\\
172 17.1666666666667\\
173 17\\
174 17\\
175 17.4\\
176 17.6\\
177 17.8\\
178 17\\
179 16\\
180 15.6\\
181 15.6\\
182 15.6\\
183 14.8\\
184 15\\
185 15.1666666666667\\
186 14.8333333333333\\
187 14.3333333333333\\
188 15.1428571428571\\
189 14.4\\
190 14.2\\
191 14.2\\
192 14.2\\
193 14.4\\
194 14\\
195 14.8\\
196 14\\
197 14.2\\
198 14.2\\
199 14.4\\
200 14.6\\
201 14.8\\
202 14.8\\
203 14.8\\
204 15\\
205 15.6\\
206 15.2\\
207 15.2\\
208 15.8\\
209 16\\
210 16\\
211 16.2\\
212 15.5\\
213 15.2\\
214 15.6\\
215 16.3333333333333\\
216 16.6666666666667\\
217 16.8\\
218 16.6\\
219 17\\
220 17.4\\
221 17.8\\
222 17.1666666666667\\
223 16\\
224 15.8\\
225 15.8\\
226 15.6\\
227 14.6\\
228 15\\
229 15.1666666666667\\
230 15\\
231 14.3333333333333\\
232 14.3333333333333\\
233 15\\
234 14.3333333333333\\
235 14.2\\
236 14.4\\
237 14.6\\
238 14.6\\
239 14.4\\
240 14\\
241 14.2\\
242 14.4\\
243 14.6\\
244 14.8\\
245 15\\
246 15.2\\
247 15.2\\
248 15.2\\
249 15.4\\
250 15.2\\
};
\addlegendentry{ILU}
\addplot [
color=mycolor2,
mark=square,
solid
]
table[row sep=crcr]{
1 10.6\\
2 11\\
3 10.4\\
4 10.4\\
5 10.5\\
6 10.5\\
7 10.4\\
8 10.4\\
9 10.4\\
10 10.2\\
11 10\\
12 10\\
13 9.8\\
14 10.2\\
15 10.4\\
16 10.4\\
17 10.6\\
18 10.6\\
19 10.6\\
20 10.4\\
21 10.4\\
22 10.4\\
23 10.4\\
24 10.4\\
25 10.4\\
26 10.4\\
27 10.8\\
28 10.8\\
29 10.6\\
30 10.6\\
31 10.6\\
32 10.6\\
33 11\\
34 11.1666666666667\\
35 11.1666666666667\\
36 11\\
37 10.8\\
38 10.6\\
39 10.5\\
40 10.1666666666667\\
41 10.4\\
42 11\\
43 11.2\\
44 11.4\\
45 11\\
46 10.8\\
47 10.4\\
48 10.6\\
49 10.5\\
50 10.5\\
51 10.6666666666667\\
52 10.4\\
53 10.4\\
54 10.4\\
55 10.2\\
56 10.2\\
57 10.2\\
58 10\\
59 10.2\\
60 10.4\\
61 10.6\\
62 10.6\\
63 10.6\\
64 10.4\\
65 10.4\\
66 10.4\\
67 10.2\\
68 10.2\\
69 10.4\\
70 10.4\\
71 10.6\\
72 10.6\\
73 10.6\\
74 10.6\\
75 10.8\\
76 11.2\\
77 11.2\\
78 11.5\\
79 11.1666666666667\\
80 11.1666666666667\\
81 10.8\\
82 11\\
83 10.6666666666667\\
84 10.5\\
85 10.6\\
86 10.6\\
87 11.4\\
88 11.6\\
89 11.5\\
90 11\\
91 10.8\\
92 10.6\\
93 10.4\\
94 10.5\\
95 10.6666666666667\\
96 10.4\\
97 10.4\\
98 10.4\\
99 10.4\\
100 10\\
101 10\\
102 10\\
103 10.2\\
104 10.4\\
105 10.4\\
106 10.6\\
107 10.6\\
108 10.4\\
109 10.4\\
110 10.4\\
111 10.4\\
112 10.4\\
113 10.4\\
114 10.6\\
115 10.8\\
116 10.8\\
117 10.8\\
118 10.8\\
119 10.8\\
120 10.6\\
121 10.8\\
122 10.8\\
123 11.1666666666667\\
124 11.1666666666667\\
125 10.8\\
126 10.8\\
127 10.5\\
128 10.3333333333333\\
129 10.4\\
130 10.6\\
131 10.8\\
132 11.4\\
133 11.5\\
134 11\\
135 10.2\\
136 10.4\\
137 10.4\\
138 10.5\\
139 10.6666666666667\\
140 10.6666666666667\\
141 10.4\\
142 10.4\\
143 10.4\\
144 10.2\\
145 10\\
146 10\\
147 10.4\\
148 10.4\\
149 10.6\\
150 10.6\\
151 10.6\\
152 10.4\\
153 10.4\\
154 10.4\\
155 10.4\\
156 10.2\\
157 10.4\\
158 10.4\\
159 10.4\\
160 10.6\\
161 10.6\\
162 10.6\\
163 10.6\\
164 11.2\\
165 11.2\\
166 11.4\\
167 11.5\\
168 11.1666666666667\\
169 10.8\\
170 10.8\\
171 10.6666666666667\\
172 10.5\\
173 10.6\\
174 10.8\\
175 11\\
176 11.4\\
177 11.6\\
178 11\\
179 11\\
180 10.4\\
181 10.4\\
182 10.5\\
183 10.5\\
184 10.4\\
185 10.4\\
186 10.4\\
187 10.2\\
188 10\\
189 10\\
190 10\\
191 10.2\\
192 10.4\\
193 10.4\\
194 10.6\\
195 10.6\\
196 10.6\\
197 10.4\\
198 10.4\\
199 10.4\\
200 10.4\\
201 10.4\\
202 10.6\\
203 10.6\\
204 10.8\\
205 10.8\\
206 10.6\\
207 10.8\\
208 10.6\\
209 10.6\\
210 11\\
211 11.1666666666667\\
212 11.1666666666667\\
213 11\\
214 10.8\\
215 10.6\\
216 10.5\\
217 10.1666666666667\\
218 10.4\\
219 11\\
220 11.2\\
221 11.6\\
222 11\\
223 10.6\\
224 10.4\\
225 10.4\\
226 10.5\\
227 10.5\\
228 10.6666666666667\\
229 10.4\\
230 10.4\\
231 10.4\\
232 10.2\\
233 10\\
234 10.2\\
235 10\\
236 10.4\\
237 10.4\\
238 10.6\\
239 10.6\\
240 10.4\\
241 10.4\\
242 10.4\\
243 10.4\\
244 10.2\\
245 10.4\\
246 10.4\\
247 10.4\\
248 10.6\\
249 10.6\\
250 10.6\\
};
\addlegendentry{Vanka};
\end{axis}
\end{tikzpicture}%
  \caption{3d fsi-3 configuration: Mean number of GMRES steps to solve
    the momentum equation \eqref{problem:1} within every Newton step
    plotted over time steps}  
  \label{fig:meangmresfsi3d}
\end{figure}
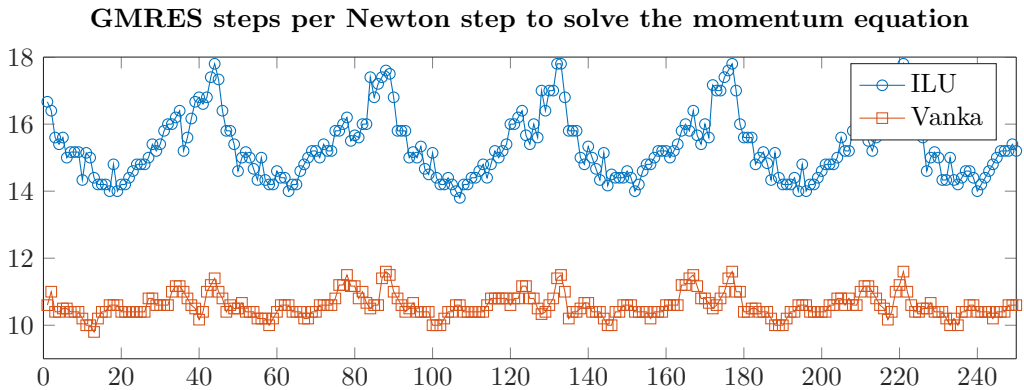

The average number of Newton steps and the average number of Jacobians
assembled for all time steps within the intervals are gathered in
Table \ref{table:Newtonsteps2d3d}. In addition, we present the mean number of
GMRES steps to solve the linearized momentum equation, problem
\eqref{problem:1}, once. The values show that the average number of
matrix assemblies in each time step can be below 1. This is due to the
approximation of the Jacobian by reassembling it, only if the
convergence rates deteriorate. Both multigrid approaches, Vanka and
ILU are very robust with regard to mesh refinement. The linear
iteration counts rise only slightly.

Figure~\ref{fig:meangmresfsi3d} shows the average number of GMRES
steps required for both Vanka and for ILU smoothing in every time
step. The values fluctuate due to the oscillatory motion of the beam.

According to Figure \ref{fig:memoryconsumtion} we need $43.88$s for
each Newton step on mesh level 6
in the 2d configuration. And according to Table
\ref{table:Newtonsteps2d3d} an average of 5.1 Newton steps. The mean
computational time per time step is  $\unit[43.88]{s} \cdot
5.1=\unit[223.49]{s}$, whereby an average of
$\unit[7.01]{s}\cdot1.23=\unit[8.66]{s}$ are used to construct the
Jacobian. Most of the computational time is spent by the linear
solver. In every Newton step the linear solver needs about
$(\unit[223.9]{s}-\unit[8.66]{s} )/{5.1}=\unit[42]{s}$. This is very
close to the value in \cite{AulisaBnaBornia2018}, where about
$\unit[46.1]{s}$ per linear solve are needed on the same level. On a
different 3d configuration with smaller deformation presented in
\cite{JodlbauerLanger}, a mesh with 
$14\cdot10^6$ degrees of freedom required $\unit[7962]{s}$ per Newton
step using a parallel block-preconditioned GMRES method on 16
cores. If we extrapolate the computational time 
in Figure~\ref{fig:memoryconsumtion}, we expect to need about $\unit[2345]{s}$
per Newton step (in single core performance). We want to highlight
that the two configurations are not directly comparable.

%

\subsubsection{Parallelization}\label{sec:parallel}

\begin{figure}[t]
  \centering
  \setlength\figureheight{0.35\textwidth} 
  \setlength\figurewidth{0.45\textwidth} 
  \small
  \begin{tikzpicture}
  \begin{groupplot}[
      group style={
        group size=1 by 2,
        xticklabels at=edge bottom,
        vertical sep=0pt
      },
      legend cell align=left,
      legend pos=north east,
      width=0.9\textwidth,
      height = 4cm,
      xmin=0, xmax= 22
    ]
    \nextgroupplot[
      ymin=200, ymax = 8000,
      ytick={2000,4000,6000,8000},
      axis x line*=top, 
      height = 3.5cm,
      axis y discontinuity=parallel,
      ]
      \addplot[blue,mark=*] plot coordinates {
        (0, 6716)
        (1, 6528)
        (2, 6557)
        (3, 5998)
        (4, 6595)
        (5, 6397)
        (6, 6283)
        (7, 5554)
        (8, 2551)
        (9, 1075)
        (10, 449)};
      \coordinate (top) at (rel axis cs:0,1);
      \addlegendentry{Patches in fluid color}
      \addplot[red,mark=triangle*] plot coordinates{(0,0)};
      \addlegendentry[red]{Patches in solid color}
    \nextgroupplot[
        ymin=0,ymax = 100,
        axis x line*=bottom,
        ytick={0,50,100},
        xlabel=Coloring,
        height = 2.5cm,
      ]
    \addplot[blue,mark=*] plot coordinates {
      (11, 63 )
      (12, 16 )
      (13, 2  ) };
    \addplot[red,mark=triangle*] plot coordinates {
      (14, 80)
      (15, 80)
      (16, 80)
      (17, 80)
      (18, 80)
      (19, 80)
      (20, 80)
      (21, 80) };
    \coordinate (bot) at (rel axis cs:1,0);
  \end{groupplot}
\end{tikzpicture}
  \caption{Coloring for avoiding memory collisions in the parallel
    Vanka smoother for the finest 3d mesh with about $3.5\cdot 10^6$
    dofs. The fluid patches consist of $108$ dofs each, while the
    solid patches couple $500$ dofs. The smallest fluid color has only
    2 patches. 
    }\label{fig:colors}
\end{figure}
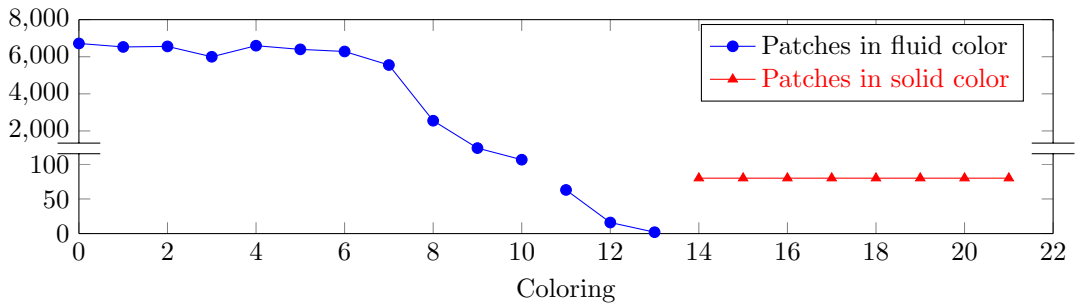

The Vanka smoother (based on a Jacobi iteration) has the advantage
that it can be easily parallelized. We introduce a cell wise coloring
of the Vanka
patches. Colors are attributed by a simple ad hoc algorithm. We run
over all patches; if a patch is not already labeled with a
color, we will label it and block all its neighbours that share a
common degree of freedom for this color. Then, we continue with the
next color. This algorithm is not optimal in terms of ``numbers of
colors'' and also not optimal in terms of ``balanced number of
elements per color'' but adequate for our purpose. 
As different patch sizes for fluid and
solid domain are used in 3d, a different color is always allocated to
fluid and solid patches, such that a good load balancing is
possible.
The finest mesh level in 3d is partitioned into 22 colors (13 within
the fluid, 9 in the solid domain),
see Figure~\ref{fig:colors}, whereby the number of patches in each color ranges
between 6\,716 and 2 within the fluid domain, and is constant with 80
pachtes per color within the solid. 
About 99\% of the fluid patches
belong to colors containing at least 500 patches, such that very little
overhead must be expected due to suboptimality of partitioning the
remaining colors (as long as a moderate number of threads is
considered). Our 
algorithms yields solid colors with 80 patches 
each. While 80 is dividable by 16, it is not dividable by 32. Hence,
the potential efficiency of functions depending on this coloring is
reduced to about 0.8 (for 32 threads).

Furthermore, we parallelized the matrix vector product. Although in
principle trivial to parallelize, we suffer
from the usually memory bandwidth restrictions that will limit
possible speedups for matrix vector products.
All parallelization is done in
\emph{OpenMP}~\cite{openmp}. We note that the parallelization is not
the focus of this work. Only first steps have been
undertaken and the implementation allows for further optimization. 

\begin{figure}[t]
  \centering
  \textbf{Mean time (in seconds) per time step}\\
  \setlength \figureheight{0.35\textwidth} 
  \setlength\figurewidth{0.43\textwidth} 
  \small
  \begin{tikzpicture}
\begin{loglogaxis}[%
width=\figurewidth,
height=\figureheight,
scale only axis,
xtick={1,2,4,8,16,32},
xticklabels={1,2,4,8,16,32},
xmin=0.75,
xmax=200,
ymin=0.5,
ymax=1000,
title style={font=\bfseries,yshift=-0.7cm,xshift=-1cm},
title={2d benchmark},
]
\addplot [
color=mycolor5,
solid,
mark=pentagon*
]
table[row sep=crcr]{
1 287.691\\2	183.46428\\4	116.22972\\8	69.26652\\16	45.6686\\32	38.01924 \\
};
\addlegendentry{all}
\addplot [
color=mycolor1,
solid,
mark=*
]
table[row sep=crcr]{
1 35.77788\\ 
2 22.17656\\	
4 13.145\\
8 9.46768\\	
16 7.66736\\
32 7.34588\\
};
\addlegendentry{NR}
\addplot [
color=mycolor2,
solid,
mark=square*
]
table[row sep=crcr]{
1 9.19512\\	2 5.05316\\	4 2.70828\\	8 1.54276\\	16 0.99924\\	32 0.8278\\
};
\addlegendentry{As}
\addplot [
color=mycolor3,
solid,
mark=diamond*
]
table[row sep=crcr]{
1 149.05136\\	2	95.44004\\	4	53.25268\\	8	29.53732\\	16	17.58708\\	32	13.91744\\	
};
\addlegendentry{MS}
\addplot [
color=mycolor4,
solid,
mark=triangle*
]
table[row sep=crcr]{
1 86.27272\\	2	52.65888\\	4	38.77592\\	8	20.44608\\	16	10.98212\\	32	7.4684\\	
};
\addlegendentry{MV}
\addplot [
color=mycolor5,
dashed,
]
table[row sep=crcr]{
1 100 \\ 32 3.125\\
};
\addlegendentry{lin}
\end{loglogaxis}
\end{tikzpicture}%
  \begin{tikzpicture}
\begin{loglogaxis}[%
width=\figurewidth,
height=\figureheight,
scale only axis,
xmin=0.75,
xmax=200,
ymin=5,
ymax=10000,
xtick={1,2,4,8,16,32},
xticklabels={1,2,4,8,16,32},
title style={font=\bfseries,yshift=-0.7cm,xshift=-1cm},
title={3d benchmark},
]
\addplot [
color=mycolor5,
solid,
mark=pentagon*
]
table[row sep=crcr]{
1 2552.34024\\
2 1435.3484\\4	813.97112\\8	556.15068\\16	328.28076\\32	245.6646\\
};
\addlegendentry{all}
\addplot [
color=mycolor1,
solid,
mark=*
]
table[row sep=crcr]{
1 285.0162\\
2 150.16108\\4	78.45576\\8	42.67888\\16	26.72912\\32	22.16548\\
};
\addlegendentry{NR}
\addplot [
color=mycolor2,
solid,
mark=square*
]
table[row sep=crcr]{
1 138.51644	\\
2 78.96392\\4	43.10224\\8	26.21932\\16	16.6624\\32	11.81628\\
};
\addlegendentry{As}
\addplot [
color=mycolor3,
solid,
mark=diamond*
]
table[row sep=crcr]{
1 998.8784	\\
2 607.98176\\4	325.57532\\8	191.25736\\16	118.6736\\32	102.13948\\
};
\addlegendentry{MS}
\addplot [
color=mycolor4,
solid,
mark=triangle*
]
table[row sep=crcr]{
1 1110.528\\
2	576.42964\\4	346.09832\\8	275.16652\\16	145.71592\\32	88.67028\\
};
\addlegendentry{MV}
\addplot [
color=mycolor5,
dashed
]
table[row sep=crcr]{
1 2000\\32 62.5\\
};
\addlegendentry{lin}
\end{loglogaxis}
\end{tikzpicture}%
  \caption{Strong scalability test. 
    Mean time per timestep (all) to compute the Newton residual
    (NR), assemble the Jacobian (As),
    multilevelsolver (MS), matrix-vector multiplication (MV) in 2d
    (left) and 3d (right) using 1-32 threads on mesh levels 6 and 3.} 
  \label{fig:time_testparallel}
\end{figure}
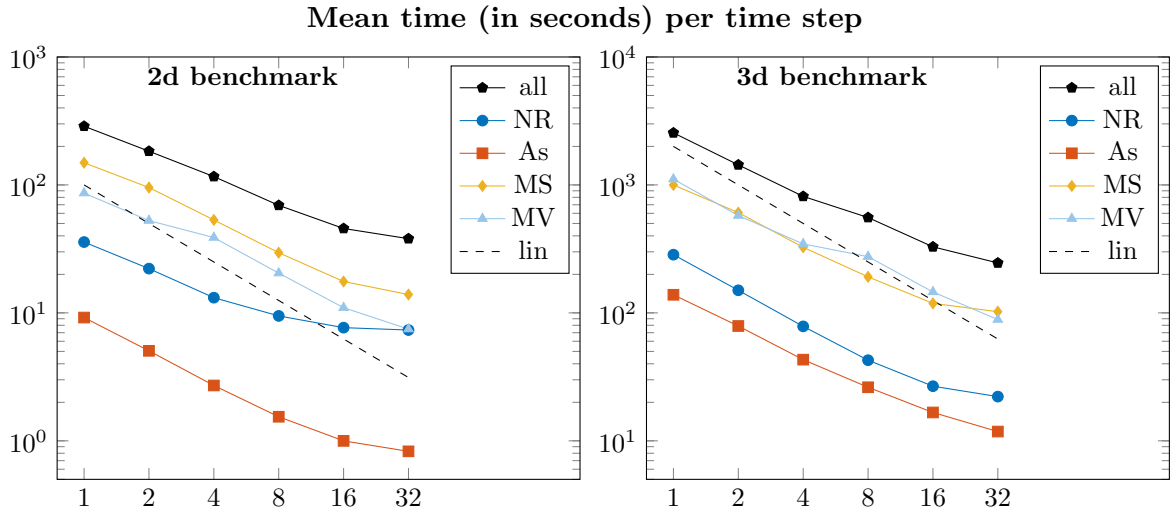

Similar to Section
\ref{sec:num-geo-mult} we recompute the 2d and 3d problem on the
time-interval $I=[9,9.5]$ with the step size  $k=0.002$s using the
finest refinement levels 6 (in 2d) and 3  (in 3d). The mean
computational time per time step on an Intel(R) Xeon(R) Gold 6150 CPU
@ 2.70GHz is given in Figure \ref{fig:time_testparallel} in a strong
scalability test. In 3d we can observe that the parallelization of all
ingredients scales rather well. If we double the number of cores the
computational time reduces by a factor of 0.57. With 32 threads we
achieve a speed up of about 10 in comparison to single core
performance. The drop in efficiency from 16 to 32 threads (in 3d) is
clearly visible in the assembly of the residual and the application of
the Vanka smoother, two functions that strongly depend on the coloring
of the patches. 

In Table~\ref{tab:parallel:dist} we show, how the distribution of the
computational time to the different ingredients develops for
an increasing numbers of threads. These results belong to the 3d
benchmark problem on the finest mesh level 3. The numbers show that
more than 80\% of the time is spend in linear algebra routines like
sparse matrix-vector products and the application of the Vanka
smoother. These operations are mainly limited by the memory
bandwidth. The very low contribution of only 5\% for the matrix
assembly could lead to the conclusion that a matrix free
implementation might be the proper choice. However, our implementation
requries less than one matrix assembly per time step in average. The
multigrid smoother however is applied many hundred times (about 5
Newton steps, 10 GMRES steps each, several Vanka steps). A matrix free
implementation on such a low number of threads would hence strongly
increase the overall time.

\begin{table}[t]
\begin{center}
  \begin{tabular}{c|ccccc}
    \toprule
    \# Threads & Total & Residual & Matrix & MV product & Vanka \\
    \midrule
    1 & 100\% & 11\%& 5\%& 44\%& 39\%\\
    4 & 100\% & 10\%& 5\%& 43\%& 40\%\\
    16 &100\% & 8 \%& 5\%& 45\%& 36\%\\
    \bottomrule
  \end{tabular}
\end{center}
\caption{Distribution of the computational time to the  main
  ingredients of the Newton-multigrid solver: integration of the
  nonlinear residual, assembly of the Jacobian, matrix-vector products
  and application of the Vanka smoother. The numbers do not add up to
  100\% as some parts, like the memory management, are not included in
  the measurement. }
\label{tab:parallel:dist}
\end{table}

\section{Summary}\label{sec:conclusion}

We have introduced a Newton multigrid framework for monolithic
fluid-structure interactions in ALE coordinates. The solver is based
on two reduction techniques in the Jacobian: first, a 
condensation of the solid deformation by representing the deformation
gradient on the velocity only and second, by skipping the ALE
derivatives within the Navier-Stokes equation. This second steps leads
to an approximated Newton method but we could show (also in
preliminary works) that the time-to-solution even benefits from this
approximation, as the computational time for assembling the ALE
derivatives is very high. The reduction has two positive effect: the
large system of 7 unknowns (in 3d) decomposes into on fluid-solid
problem in pressure and velocity with 4 unknowns and two partitioned
systems with 3 unknowns each for solving solid and fluid
deformation. The second effect is the better conditioning of the
coupled system that allows for the use of very simple multigrid
smoothers that are easy to parallelize. Also, while ILU smoothing
applied to the monolithic system was
not convergent in our previous contribution~\cite{Richter2015}, is
performed well for smoothing the global momentum equations. 
Combined with first steps of parallelization and in comparison to our
past approaches based on a monolithic solution of the complete
pressure-velocity-deformation system and partitioned smoothers and
also in comparison to approaches presented in literature we could
significantly reduce the computational time.

As basis for future benchmarking of 3d fluid-structure interactions we
presented an extension of the 2d benchmark problems by Hron and
Turek~\cite{HronTurek2006} that is by far more challenging (due to
larger deformations and a strong dynamic behavior) as compared to a
first test case introduced in our past work~\cite{Richter2015} which
has also been considered in~\cite{AulisaBnaBornia2018,JodlbauerLanger}
in very similar studies. It 
will still require further effort to establish reference values for
this new 3d benchmark case. 

Our work includes some first simple steps of parallelization which
have to be extended in future work. In particular, in order to
overcome the memory bandwidth limitations which are common in such memory
extensive computations, distributed memory paradigms have to be
incorporated~\cite{KimmritzRichter2010}. Further, some benefit can be
expected by using GPU acceleration for matrix vector product and Vanka
smoother.

\section*{Acknowledgements}
Both authors acknowledge the financial support by the Federal Ministry of
Education and Research of Germany, grant number 05M16NMA, TR
acknowledges the support of the GRK 2297 MathCoRe, funded by the Deutsche
Forschungsgemeinschaft, grant number 314838170.


\begin{thebibliography}{10}
\providecommand{\url}[1]{{#1}}
\providecommand{\urlprefix}{URL }
\expandafter\ifx\csname urlstyle\endcsname\relax
  \providecommand{\doi}[1]{DOI~\discretionary{}{}{}#1}\else
  \providecommand{\doi}{DOI~\discretionary{}{}{}\begingroup
  \urlstyle{rm}\Url}\fi

\bibitem{MUMPS:2}
Amestoy, P.R., Guermouche, A., L'Excellent, J.Y., Pralet, S.: Hybrid scheduling
  for the parallel solution of linear systems.
\newblock Parallel Computing \textbf{32}(2), 136--156 (2006)

\bibitem{AulisaBnaBornia2018}
Aulisa, E., Bna, S., Bornia, G.: A monolithic ale newton-krylov solver with
  multigrid-richardson-schwarz preconditioning for incompressible
  fluid-structure interaction.
\newblock Computers \& Fluids \textbf{174}, 213--228 (2018)

\bibitem{BeckerBraack2000a}
Becker, R., Braack, M.: Multigrid techniques for finite elements on locally
  refined meshes.
\newblock Numerical Linear Algebra with Applications \textbf{7}, 363--379
  (2000).
\newblock Special Issue

\bibitem{BeckerBraack2001}
Becker, R., Braack, M.: A finite element pressure gradient stabilization for
  the {S}tokes equations based on local projections.
\newblock Calcolo \textbf{38}(4), 173--199 (2001)

\bibitem{Gascoigne3D}
Becker, R., Braack, M., Meidner, D., Richter, T., Vexler, B.: The finite
  element toolkit \textsc{Gascoigne}.
\newblock \textsc{http://www.gascoigne.uni-hd.de}

\bibitem{BraackRichter2006d}
Braack, M., Richter, T.: Solutions of {3D} {N}avier-{S}tokes benchmark problems
  with adaptive finite elements.
\newblock Computers and Fluids \textbf{35}(4), 372--392 (2006)

\bibitem{BrummelenZeeBorst2008}
Brummelen, E., Zee, K., Borst, R.: Space/time multigrid for a
  fluid-structure-interaction problem.
\newblock Applied Numerical Mathematics \textbf{58}(12), 1951--1971 (2008)

\bibitem{BungartzSchaefer2006}
Bungartz, H.J., Sch\"afer, M. (eds.): Fluid-Structure Interaction. Modelling,
  Simulation, Optimisation, \emph{Lecture Notes in Computational Science and
  Engineering}, vol.~53.
\newblock Springer (2006).
\newblock ISBN-10: 3-540-34595-7

\bibitem{BungartzSchaefer2010}
Bungartz, H.J., Sch\"afer, M. (eds.): Fluid-Structure Interaction II.
  Modelling, Simulation, Optimisation.
\newblock Lecture Notes in Computational Science and Engineering. Springer
  (2010)

\bibitem{CausinGerbeauNobile2005}
Causin, P., Gereau, J., Nobile, F.: Added-mass effect in the design of
  partitioned algorithms for fluid-structure problems.
\newblock Comput. Methods Appl. Mech. Engrg. \textbf{194}, 4506--4527 (2005)

\bibitem{CrosettoDeparis}
Crosetto, P., Deparis, S., Fourestey, G., Quarteroni, A.: Parallel algorithms
  for fluid-structure interaction problems in haemodynamics.
\newblock SIAM Journal on Scientific Computing \textbf{33}(4), 1598--1622
  (2011).
\newblock \doi{10.1137/090772836}

\bibitem{Davis2014}
Davis, T.: Umfpack, an unsymmetric-pattern multifrontal method.
\newblock ACM Transactions on Math. Soft. \textbf{30}(2), 196--199 (2014)

\bibitem{DeparisForti}
Deparis, S., Forti, D., Grandperrin, G., Quarteroni, A.: Facsi: A block
  parallel preconditioner for fluid-structure interaction in hemodynamics.
\newblock Journal of Computational Physics \textbf{327}, 700 -- 718 (2016).
\newblock \doi{https://doi.org/10.1016/j.jcp.2016.10.005}.
\newblock
  \urlprefix\url{http://www.sciencedirect.com/science/article/pii/S0021999116304983}

\bibitem{Failer2017}
Failer, L.: Optimal control for time dependent nonlinear fluid-structure
  interaction.
\newblock Ph.D. thesis, Technische Universit\"at M\"unchen (2017)

\bibitem{FailerWick}
Failer, L., Wick, T.: Adaptive time-step control for nonlinear fluid?structure
  interaction.
\newblock Journal of Computational Physics \textbf{366}, 448 -- 477 (2018)

\bibitem{FernandezGerbeau2009}
Fern{\'a}ndez, M., Gerbeau, J.F.: Algorithms for fluid-structure interaction
  problems.
\newblock In: L.~Formaggia, A.~Quarteroni, A.~Veneziani (eds.) Cardiovascular
  Mathematics: Modeling and simulation of the circulatory system, \emph{MS \&
  A}, vol.~1, pp. 307--346. Springer (2009)

\bibitem{FernandezMoubachir2005}
Fern{\'a}ndez, M., Moubachir, M.: A newton method using exact jacobians for
  solving fluid-structure coupling.
\newblock Computers and Structures \textbf{83}, 127--142 (2005)

\bibitem{Frei2016}
Frei, S.: Eulerian finite element methods for interface problems and
  fluid-structure interactions.
\newblock Ph.D. thesis, Universit\"at Heidelberg (2016).
\newblock Doi:10.11588/heidok.00021590

\bibitem{GeeKuettlerWall2010}
Gee, M., K{\"u}ttler, U., Wall, W.: Truly monolithic algebraic multigrid for
  fluid-structure interaction.
\newblock Int. J. Numer. Meth. Engrg. \textbf{85}, 987--1016 (2010)

\bibitem{eigenweb}
Guennebaud, G., Jacob, B., et~al.: Eigen v3.
\newblock http://eigen.tuxfamily.org (2010)

\bibitem{HeilHazelBoyle2008}
Heil, M., Hazel, A., Boyle, J.: Solvers for large-displacement fluid-structure
  interaction problems: Segregated vs. monolithic approaches.
\newblock Computational Mechanics \textbf{43}, 91--101 (2008)

\bibitem{HeywoodRannacherTurek1992}
Heywood, J., Rannacher, R., Turek, S.: Artificial boundaries and flux and
  pressure conditions for the incompressible {N}avier-{S}tokes equations.
\newblock Int. J. Numer. Math. Fluids. \textbf{22}, 325--352 (1992)

\bibitem{HronTurek2006a}
Hron, J., Turek, S.: A monolithic {FEM}/{M}ultigrid solver for an {ALE}
  formulation of fluid-structure interaction with applications in biomechanics.
\newblock In: H.J. Bungartz, M.~Sch{\"a}fer (eds.) Fluid-Structure Interaction:
  Modeling, Simulation, Optimization, Lecture Notes in Computational Science
  and Engineering, pp. 146--170. Springer (2006)

\bibitem{HronTurek2006}
Hron, J., Turek, S.: Proposal for numerical benchmarking of fluid-structure
  interaction between an elastic object and laminar incompressible flow.
\newblock In: H.J. Bungartz, M.~Sch{\"a}fer (eds.) Fluid-Structure Interaction:
  Modeling, Simulation, Optimization, Lecture Notes in Computational Science
  and Engineering, pp. 371--385. Springer (2006)

\bibitem{JodlbauerLanger}
Jodlbauer, D., Langer, U., Wick, T.: Parallel block-preconditioned monolithic
  solvers for fluid-structure interaction problems.
\newblock International Journal for Numerical Methods in Engineering
  \textbf{117}(6), 623--643 (2019)

\bibitem{KimmritzRichter2010}
Kimmritz, M., Richter, T.: Parallel multigrid method for finite element
  simulations of complex flow problems on locally refined meshes.
\newblock Numerical Linear Algebra with Applications \textbf{18}(4), 615--636
  (2010)

\bibitem{openmp}
Klemm, M., Supinski, B. (eds.): OpenMP Application Programming Interface
  Specification Version 5.0.
\newblock Independently published (2019)

\bibitem{LangerYang2017}
Langer, U., Yang, H.: Recent development of robust monolithic fluid-structure
  interaction solvers.
\newblock In: Fluid-Structure Interactions. Modeling, Adaptive Discretization
  and Solvers, \emph{Radon Series on Computational and Applied Mathematics},
  vol.~20. de Gruyter (2017)

\bibitem{Molnar2015}
Molnar, M.: {S}tabilisierte {F}inite {E}lemente f\"ur {S}tr\"omungsprobleme auf
  bewegten {G}ebieten.
\newblock Master's thesis, Universit\"at Heidelberg (2015)

\bibitem{Pironneau2016}
Pironneau, O.: An energy preserving monolithic eulerian fluid-structure
  numerical scheme.
\newblock Chinese Annals of Mathematics \textbf{39} (2016).
\newblock Preprint at arXiv:1607.08083

\bibitem{Pironneau2019}
Pironneau, O.: An Energy stable Monolithic {E}ulerian Fluid-Structure Numerical
  Scheme with compressible materials (2019).
\newblock Https://arxiv.org/abs/1607.08083

\bibitem{Richter2015}
Richter, T.: A monolithic geometric multigrid solver for fluid-structure
  interactions in {ALE} formulation.
\newblock Int. J. Numer. Meth. Engrg. \textbf{104}(5), 372--390 (2015)

\bibitem{Richter2017}
Richter, T.: Fluid-structure Interactions. Models, Analysis and Finite
  Elements, \emph{Lecture Notes in Computational Science and Engineering}, vol.
  118.
\newblock Springer (2017)

\bibitem{RichterWick2010}
Richter, T., Wick, T.: Finite elements for fluid-structure interaction in {ALE}
  and {F}ully {E}ulerian coordinates.
\newblock Comput. Methods Appl. Mech. Engrg. \textbf{199}(41-44), 2633--2642
  (2010)

\bibitem{RichterWick2015_time}
Richter, T., Wick, T.: On time discretizations of fluid-structure interactions.
\newblock In: T.~Carraro, M.~Geiger, S.~K\"orkel, R.~Rannacher (eds.) Multiple
  Shooting and Time Domain Decomposition Methods, \emph{Contributions in
  Mathematical and Computational Science}, vol.~9, pp. 377--400. Springer
  (2015)

\bibitem{TurekHronMadlikRazzaqWobkerAcker2010}
Turek, S., Hron, J., Madlik, M., Razzaq, M., Wobker, H., Acker, J.: Numerical
  simulation and benchmarking of a monolithic multigrid solver for
  fluid--structure interaction problems with application to hemodynamics.
\newblock Tech. rep., Fakult\"{a}t f\"{u}r Mathematik, TU Dortmund (2010).
\newblock Ergebnisberichte des Instituts f\"{u}r Angewandte Mathematik, Nummer
  403

\bibitem{TurekRivkindHronGlowinski2006}
Turek, S., Rivkind, L., Hron, J., Glowinski, R.: Numerical study of a modified
  time--stepping theta--scheme for incompressible flow simulations.
\newblock Journal of Scientific Computing \textbf{28}(2--3), 533--547 (2006)

\bibitem{Wall1999}
Wall, W.: Fluid-structure interaction with stabilized finite elements.
\newblock Ph.D. thesis, University of Stuttgart (1999).
\newblock Urn:nbn:de:bsz:93-opus-6234

\bibitem{YirgitSchaeferHeck2008}
Yirgit, S., Sch\"afer, M., Heck, M.: Grid movement techniques and their
  influence on laminar fluid-structure interaction rpoblems.
\newblock J. Fluids and Structures \textbf{24}(6), 819--832 (2008)

\bibitem{ZeeBrummelenBorst2010a}
Zee, K., Brummelen, E., Borst, R.: Goal-oriented error estimation and
  adaptivity for free-boundary problems: The domain-map linearization approach.
\newblock SIAM J. on Scientific Computing \textbf{32}(2), 1074 -- 1092 (2010)

\bibitem{ZeeBrummelenBorst2010b}
Zee, K., Brummelen, E., Borst, R.: Goal-oriented error estimation and
  adaptivity for free-boundary problems: The shape-linearization approach.
\newblock SIAM J. on Scientific Computing \textbf{32}(2), 1093--1118 (2010)

\end{thebibliography}
\end{document}